\documentclass[moor]{informs1}      
\usepackage[sort]{cite}
\usepackage{bm}
\usepackage{algorithmic} 
\usepackage{algorithm} 
\usepackage{titletoc}
\usepackage{graphics}
\usepackage{graphicx}
\usepackage{tikz}
\usepackage{framed}
\usepackage{bm}
\usepackage{pgfplots}
\usepackage{booktabs}
\usepackage{mathtools}
\usepackage{siunitx}
\usepackage{changepage}
\usepackage{placeins}
\usepackage{caption}
\usepackage{subcaption}
\usepackage{multicol}
\usepackage{setspace}
\usepackage[shortlabels]{enumitem}
\setlist[enumerate]{{\rm(a)}}

\usepackage{tcolorbox}     
\tcbset{%
sharp corners,
boxrule=0.5pt,
before={\par\pagebreak[0]\medskip\parindent=0pt},
after={\noindent}
}

\def\argmin{ \mathop{{\rm argmin}}}

\newcommand{\co}{\mathrm{conv}\,}
\newcommand{\cl}{\mathrm{cl}\,}

\newcommand{\dom}{\mathrm{dom}\,}
\newcommand{\epi}{\mathrm{epi}\,}

\newcommand{\gph}{\mathrm{gph}\,}
\newcommand{\id}{\mathrm{id}\,}
\newcommand{\inter}{\mathrm{int}\,}

\newcommand{\ri}{\mathrm{ri}\,}

\newcommand{\sgn}{\mathrm{sgn}}

\def\Limsup{\mathop{{\rm Lim}\,{\rm sup}}}

\newcommand{\p}{\partial}

\DeclareMathOperator{\rank}{rank}
\DeclareMathOperator{\range}{range}
\newcommand{\lev}[2]{\mathrm{lev}_{#1}#2}

\newcommand{\R}{\mathbb{R}}

\newcommand{\pLwf}{p_{L,\omega,f}}
\newcommand{\pwf}{p_{\omega,f}}

\newcommand{\Pwf}{P_{\omega,f}}

\def\pto {\overset{p}{\to}} 
\def\eto {\overset{e}{\to}} 
\def\cto {\overset{c}{\to}}

\newcommand{\rbar}{\overline{\mathbb R}}

\newcommand{\bR}{\mathbb{R}}

\newcommand{\IFF}{\quad\Longleftrightarrow\quad}

\newcommand{\bB}{\mathbb{B}}
\newcommand{\bN}{\mathbb{N}}
\newcommand{\bE}{\mathbb{E }}

\newcommand{\ip}[2]{\left\langle #1,\, #2\right\rangle}
\newcommand{\half}{\frac{1}{2}}

\newcommand{\set}[2]{\left\{#1\,\left\vert\; #2\right.\right\}}

\newcommand{\sig}{\sigma}

\newcommand{\lam}{\lambda}

\newcommand{\cD}{\mathcal{D}}

\newcommand{\UND}{\ \mbox{ and }\ }

\newcommand{\defterm}[1]{\emph{#1}}
\DeclareMathOperator{\prox}{\mathit{P}\!}
\DeclareMathOperator{\proj}{\mathit{P}\!}

\usepackage{natbib}
 \NatBibNumeric
 \def\bibfont{\small}%
 \bibpunct[, ]{[}{]}{,}{n}{}{,}%

\usepackage[colorlinks=true,breaklinks=true,bookmarks=true,urlcolor=blue,
     citecolor=blue,linkcolor=blue,bookmarksopen=false,draft=false]{hyperref}

\def\EMAIL#1{\href{mailto:#1}{#1}}
\def\URL#1{\href{#1}{#1}}         

\usepackage[capitalise,nameinlink]{cleveref}
\crefformat{equation}{(#2#1#3)}

\TheoremsNumberedThrough     

\EquationsNumberedThrough    

\renewcommand{\theARTICLETOP}{} 

\begin{document}


 \RUNAUTHOR{M.P.~Friedlander, A.~Goodwin, and T.~Hoheisel}

\RUNTITLE{From perspective maps to epigraphical projections}

 \TITLE{From perspective maps to epigraphical projections}

\ARTICLEAUTHORS{
\AUTHOR{Michael P. Friedlander}
\AFF{Department of Computer Science/Department of Mathematics,
University of British Columbia\\
2366 Main Mall
Vancouver, BC, V6T 1Z4, Canada\\
 \EMAIL{michael@friedlander.io}, \URL{https://friedlander.io}}

\AUTHOR{Ariel Goodwin}
\AFF{Department of Mathematics and Statistics, McGill University\\ 805 Sherbrooke St West, Montr\'eal, Qu\'ebec,  H3A 0B9, Canada\\
\EMAIL{ariel.goodwin@mail.mcgill.ca}, \URL{https://github.com/arielgoodwin}}

\AUTHOR{Tim Hoheisel}
\AFF{Department of Mathematics and Statistics, McGill University\\805 Sherbrooke St West, Montr\'eal, Qu\'ebec,  H3A 0B9, Canada\\
\EMAIL{tim.hoheisel@mcgill.ca}, \URL{https://www.math.mcgill.ca/hoheisel/}}


\bigskip

\it Dedicated to  James V. Burke, our  collaborator and friend,  on the occasion of his 65th birthday
} 

\ABSTRACT{%
The projection onto the epigraph or a level set of a closed  proper convex function can be achieved by finding a root of  a scalar equation  that involves the proximal operator as a function of the proximal parameter.
This paper develops the variational analysis of this scalar equation. The approach is based on a study of the variational-analytic properties of general convex optimization problems that are (partial) infimal projections of the the sum of the function in question and the perspective map of a convex kernel. When the kernel is the Euclidean norm squared, the solution map corresponds to the proximal map, and thus the variational properties derived for the general case apply to the proximal case. Properties of the value function and the corresponding solution map---including local Lipschitz continuity, directional differentiability, and semismoothness---are derived. An SC$^1$ optimization framework for computing epigraphical and level-set projections is thus established. Numerical experiments on 1-norm projection illustrate the effectiveness of the approach as compared with specialized algorithms.%
}

\MSCCLASS{52A4, 65K10, 90C25, 90C46}

\maketitle

\section{Introduction}

\begin{figure}[t]
  \centering
  \includegraphics{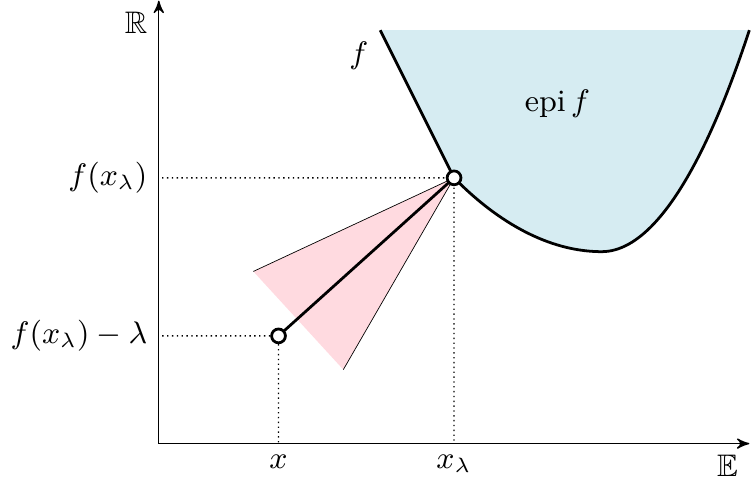}
  \caption{The proximal map $x_\lambda:=P_\lambda f(x)$ corresponds to the projection of the pair $(x,\,f(x_\lambda)-\lambda))$ onto the epigraph of $f$; see~\eqref{eq:EpiProx}.\label{fig:epi-projection}}
\end{figure}

The Moreau proximal map of a closed proper convex function $f$ that maps a finite-dimensional Euclidean space $\bE_f$ to $\rbar:=\R\cup\{+\infty\}$ is given by the minimizing set 
\begin{equation*}
P_{\lambda} f( x)=\argmin_{u\in \bE_f}\left\{ f(u)+(1/2\lambda) \|x-u\|^2\right\}\quad (\lambda>0).
\end{equation*}
The proximal map is a central operation of algorithms for nonsmooth optimization, including first-order methods such as proximal gradient and operator splitting~\cite{parikh2013proximal,BaC 17}. Geometrically, the proximal map corresponds to the Euclidean  projection  $\proj_{\epi f}$ onto the epigraph $\epi f$; see \cref{fig:epi-projection}. Indeed, for all positive $\lambda$ and  $x_\lambda:=\prox_{\lambda}f(x)$, 
\begin{equation}\label{eq:EpiProx}
  \big(x_\lambda,\, f(x_\lambda)\big)
  = P_{\epi f}(x,\,f(x_\lambda)-\lambda).
\end{equation}
Thus, the projection of an arbitrary point $(x,\alpha)\in\bE_f\times\R\not\in\epi f$ corresponds to the proximal map of the base point $x$ using the parameter $\lambda$ that is the unique positive root of the function
\begin{equation}\label{eq:root}
0<\lambda\mapsto f(x_\lambda)-\lambda-\alpha.
\end{equation}
This connection between epigraphical projection and the proximal map---described by Beck~\cite{Bec 17}, Bauschke and Combettes~\cite[Section 29.5]{BaC 17}, Chierchia et al.~\cite[Proposition~1]{CPP 15}, and Meng et al.~\cite{MSZ 05, MZGS 08}---is a defining feature of a class of \emph{epigraphical first-order methods} for structured convex optimization over $\bE_f$ that operate through a sequence of projections onto the epigraphs of the underlying functions. In effect, these methods operate on an equivalent optimization problem over $\bE_f\times\R$~\cite{CPP 15, WWK 16,TKC 14, TBKC 14}.

This paper develops a general analysis that provides, among other things, the variational properties of the maps
\begin{align*}
(x,\,\lambda)\mapsto x_\lambda:= P_\lambda f(x)
\quad\UND\quad
(x,\,\lambda)\mapsto f(x_\lambda),
\end{align*}
defined on $\bE_f\times\bR$. This analysis and its supporting calculus allows us to determine the sensitivity of the epigraphical projection with respect to the simultaneous variation of the base point $x$ and the scaling parameter $\lambda$.  Although the resulting mathematical statements are key for our deeper understanding of epigraphical first-order methods, the overall analysis applies much more generally.

The approach we take is based on the variational analysis of the optimal value function
\begin{tcolorbox}
\begin{equation}\label{eq:PerOpt}
\pLwf:(x,\,\lambda)\in\bE_x\times\R\mapsto\inf_{u\in\bE_f}
	 f(u)+\omega^\pi(L(u,x),\,\lambda)
\end{equation}	
\end{tcolorbox}
and its corresponding solution map. Here, $L$ is a linear map, and the perspective transform $\omega^\pi$ of a closed proper convex function $\omega$ is defined by $\epi\omega^\pi=\cl\R_+(\epi\omega\times\{1\})$. When the linear map $L$ is defined as $(u,\,x)\mapsto x-u$, the value function~\eqref{eq:PerOpt} is the infimal convolution of the functions $f$ and $\omega^\pi(\cdot,\,\lambda)$. For this reason, we refer to this value function as the \emph{generalized convolution} of these two functions.

The convex calculus we establish in \cref{sec:Perspective} for the analysis of the generalized convolution~\eqref{eq:PerOpt} provides a key tool for understanding several important cases. These include the variational properties of infimal convolution (\cref{sec:InfConv}); parametric constrained optimization (\cref{sec:relax}); the Moreau envelope of a convex function and the corresponding proximal map (\cref{sec:prox}); and epigraphical and level-set projections, including an {\em SC${^1}$ optimization}~\cite{FaP 03, PaQ 95} method for numerically evaluating these projections (\cref{sec:EpiLevel}).

\subsection{Contributions and related work}   

The perspective map used in generalized convolution~\eqref{eq:PerOpt} first appears in Rockafellar~\cite[Corollary 13.5.1]{Roc 70}, without a particular name attached to it. More recently, Combettes~\cite{Com 18}, Combettes and M\"uller~\cite{CoM 18, CoM 20}, and Aravkin et al.~\cite{ABD 18}, describe in detail the properties and applications of this map. Our systematic study of parametric optimization problems with perspective maps, outlined in~\cref{sec:Perspective}, appears to be new.

\subsubsection{Infimal convolution} \Cref{sec:InfConv} establishes the variational properties of infimal convolution, which occurs when the map $L$ is $(u,x)\mapsto x-u$. These results complement the functional smoothing framework described by Beck and Teboulle~\cite[Section 4.1]{BeT 12} and Burke and Hoheisel~\cite{BuH 13, BuH 17}, wherein a smooth approximation to a function $f$ is constructed through the infimal convolution with the perspective map of a smooth and strongly convex regularizer $\omega$. Bougeard et al.~\cite{BPP 91} and Str\"omberg~\cite{Str 96} provide early contributions to this topic. \cref{th:InfConvSol} describes the Lipschitzian properties of the corresponding optimal solution map---as a function of $(x,\lambda)$.   \cref{cor:SemismoothP} establishes sufficient conditions for this solution map to be {\em semismooth*}~\cite{GfO 19}. These conditions hold, for instance, when $f$ is piecewise linear-quadratic. This analysis complements the study of the proximal case by Meng et al.~\cite{MSZ 05, MZGS 08} and Milzarek~\cite{Mil 16}.

\subsubsection{Parametric constrained optimization} A general form of parametric constrained optimization occurs when we specialize the convolution kernel $\omega$ in~\eqref{eq:PerOpt} to be the indicator function to a closed convex set. \Cref{sec:relax} focuses the variational analysis of the generalized convolution operation to obtain formulas for the sensitivity of the optimal value of parametric optimization problems with relaxed linear constraints. This analysis includes perturbations to the relaxation parameter and to the right-hand side.

\subsubsection{Moreau envelope and proximal map} In \cref{sec:prox} we further focus our analysis of infimal convolution on the \emph{proximal case}, which occurs when $\omega=\frac12\|\cdot\|_2^2$. Here we develop the variational properties of the Moreau envelope and the associated proximal map as a function of the base point $x$ and the proximal parameter $\lambda$, simultaneously.  We also establish conditions under which the proximal map is semismooth*. Special attention is given to the limiting properties as $\lambda\downarrow 0$ (\cref{prop:VarConv,prop:ProxConv}) and to  continuity and smoothness properties of the proximal map (\cref{cor:Lip,cor:DirDiff,Prop:Semi*P}). Milzarek's dissertation \cite{Mil 16} includes a related analysis that generalizes the proximal parameter $\lambda$ to a positive-definite matrix, but makes no statements regarding the limiting case where $\lambda$ (or its matrix counterpart) vanishes, as we do in our general analysis. See also Attouch's seminal monograph~\cite{Att 84}. 

\subsubsection{Proximal value map}

In \cref{sec:eta} we describe the main continuity properties of the \emph{proximal value function}
\begin{equation}\label{eq:proximal-value-fn}
	0<\lambda\mapsto f(P_\lambda f(\bar x)),
\end{equation}
where $\bar x\in\bE$ is held fixed. \Cref{cor:eta} establishes its Lipschitzian properties and \cref{cor:Phi} characterizes it as the derivative of the map $\lambda\mapsto \lambda e_\lambda f(\bar x)$ on $\R_{++}$.  \Cref{prop:SemiEta} describes sufficient conditions under which the proximal value function is semismooth.

\subsubsection{Post compositions, and  epigraphical and level-set projection}

We use our analysis of the proximal value function~\eqref{eq:proximal-value-fn} to establish, via \cref{prop:PostComp}, novel variational formulas for the Moreau envelope and proximal map of {\em post-compositions}, i.e., functions of the form $g\circ \psi$, where the scalar function $g$ is increasing and convex, and $\psi$ is  closed proper convex. As a consequence, \cref{cor:Epi} provides a refined version of the epigraphical projection conditions in~\eqref{eq:EpiProx}, including analogous results for the projection onto the level set of $f$ (\cref{cor:LevelSet}). This analysis does not require the function to be finite-valued, and extends existing results~\cite{BaC 17,Bec 17}. Importantly, \cref{cor:Epi} shows that the root of the aligning equation~\eqref{eq:root} coincides with the unique minimizer of a strongly convex scalar optimization problem. It follows from \cref{prop:SemiEta} that the objective for this problem is continuously differentiable with a locally Lipschitz derivative.  We use this latter property to derive a novel SC$^1$ optimization method to find the root of the function~\eqref{eq:root} and its analog in the level-set case. Numerical experiments in \cref{sec:numerics} show that for projection onto the 1-norm unit ball, the resulting SC$^1$ method is competitive with two specialized state-of-the-art methods: CONDAT~\cite{Con 16} and IBIS~\cite{Liu 09}.

\subsection{Notation}

Let $\Gamma_0(\bE)$ denote the set of functions $f:\bE\to\rbar$ that are proper closed convex, i.e., the epigraph $\epi f=\set{(x,\alpha)\in\bE\times\R}{f(x)\le\alpha}$ contains no vertical lines and is closed convex. Its {\em level sets}  are given by $\lev{\alpha}{f}:=\set{x\in\bE}{f(x)\leq \alpha}$. The Fenchel conjugate of any function $f:\bE\to\rbar$ is $f^*(y) = \sup_{x\in\bE}\{\ip y x - f(x)\}$. The Jacobian of a differentiable map $F:\R^n\to \R^m$ at $x\in \R^n$ is denoted by $F'(x)$.   
We denote the Euclidean projection of $\bar x$ onto $C$ by $P_C(\bar x)$. Throughout, fractions such as $(1/(2\lambda))$ are abbreviated as $(1/2\lambda)$.

For a set $C\subset\bE$,  its {\em indicator function} is $\delta_C:\bE\to\rbar$ given by $ \delta_C(x):=0$ if $x\in C$ and $\delta_C(x)=+\infty$ otherwise. The subdifferential of $\delta_C$ is the {\em normal cone} of $C$, i.e., $N_C(\bar x):=\p \delta_C(\bar x):=\set{v\in\bE}{\ip{v}{x-\bar x}\leq 0\ (x\in C)}$, which is empty if $\bar x\not\in C$. The relative interior of $C$ is the set $\ri C$~\cite[Section~6]{Roc 70}, and the {\em horizon cone}  is $C^\infty:=\set{v\in \bE}{\exists \{\lambda_k\}\downarrow 0,\ \{x_k\in C\}: \lambda_kx_k\to v}$. The {\em horizon function} of $f\in \Gamma_0(\bE)$ is the  closed proper convex and positively homogeneous function $f^\infty:\bE\to\rbar$ defined via $\epi f^\infty=(\epi f)^\infty$.


Let  $f_k:\bE\to \rbar$. Then we say that the sequence $\{f_k\}$ {\em epi-converges} to a function $f:\bE\to\rbar$ if 
\begin{equation*}\label{eq:epi1}
\forall x\in\bE:\,\left\{\begin{array}{ll}  \forall \{x_k\}\to x: &    \liminf_{k\to\infty} f_k(x_k) \geq f(x),\\ \exists \{x_k\}\to x: &   \limsup_{k\to\infty} f_k(x_k)\leq f(x), \end{array}\right.
\end{equation*}
and we write $f_k\eto f$.
%
The sequence $\{f_k\}$ is said to converge {\em continuously} to $f$ if 
\[
\lim_{k\to\infty} f_k(x_k)=f(x) \quad
\forall x\in\bE\ \mbox{and}\ \{x_k\}\to x,
\]
and we write $f_k\cto f$.
Furthermore, $\{f_k\}$ is said to converge {\em pointwise} to $f$ if 
$$
\lim_{k\to\infty} f_k(x)=f(x)\quad\forall x\in\bE,
$$
and we write $f_k\pto f$.
We extend these notions  to families of functions 
$\{f_\lambda\}_{\{\lambda\downarrow 0\}}$ via
\[
f_\lambda \overset{\xi}{\to} f\quad  :\Longleftrightarrow \quad\forall\{\lambda_k\}\downarrow 0:\; f_{\lambda_k} \overset{\xi}{\to} f\quad (\xi\in \{p,e,c\}).
\]

\section{Properties of the perspective map}

The perspective map $\omega^\pi$ that appears in the generalized infimal convolution \cref{eq:PerOpt} provides a mechanism for controlling, through the parameter $\lambda$, the degree to which the functions $f$ and $\omega$ are combined. Beck and Teboulle~\cite{BeT 12} and Burke and Hoheisel~\cite{BuH 13} promoted this technique for generating smooth approximations to nonsmooth functions.

We work with the following definition of the perspective map of $\omega$, which appears in Rockafellar~\cite[Corollary~13.5.1]{Roc 70}:
\begin{equation}\label{eq:OmegaPersp} 
\omega^\pi:(z,\lambda)\in\bE_\omega\times\R\mapsto
\begin{cases}
	\lambda \omega\left(z/\lambda\right) & \mbox{if $\lambda >0$,}\\
	\omega^\infty(z) & \mbox{if $\lambda=0$,}\\
    +\infty & \mbox{if $\lambda<0$.}
\end{cases}
\end{equation}
For positive values of the parameter $\lambda$, the perspective map corresponds to \defterm{epi-multiplication}:
\[
(\lambda \star \omega)(x):=\lambda\omega\left(x/\lambda\right).
\]
The following result confirms the consistency of the perspective map~\eqref{eq:OmegaPersp} as the parameter $\lambda$ decreases towards zero.

\begin{lemma}[Variational convergence of epi-multiplication]\label{lem:EpiMultConv} Let $\phi\in \Gamma_0(\bE)$. Then as $\lambda\downarrow0$, $(\lambda\star\phi)(x) \to \phi^\infty(x)$ for all $x\in\dom\phi$, and $\lambda\star\phi \eto \phi^\infty$.
\end{lemma}

\proof{Proof.}
The pointwise convergence of $(\lambda\star\phi)$ over $\dom\phi$ follows from~\cite[Corollary~8.5.2]{Roc 70}. To prove epi-convergence, observe that, for all $\lambda>0$ and $x\in\bE$,   
 $
 (\lambda\star \phi)(x)=\phi^\pi(x,\lambda).
 $
Hence,
\[
\liminf_{\overset{x\to \bar x}{\lambda\downarrow 0}}\, (\lambda\star \phi)(x)=\liminf_{\overset{x\to \bar x}{\lambda\downarrow 0}} \phi^\pi(x,\lambda)\geq  \phi^\pi(\bar x,0) =  \phi^\infty(\bar x)
\quad \forall \bar x\in\bE,
\]
where the inequality follows because $\omega^\pi$ is a support function~\cite[Corollary~13.5.1]{Roc 70} and thus closed~\cite[Proposition~2.1.2]{HUL 01}.

Fix any sequence $\{\lambda_k\}\downarrow 0$ and take $\bar x\in \dom \phi$. Then $(\lambda_k\star \phi)(\bar x)\to \phi^\infty(\bar x)$. Hence, in particular, with $x_k:=\bar x\;(k\in \mathbb N)$,
\begin{equation}\label{eq:limsup}
\limsup_{k\to\infty}\, (\lambda_k\star \phi)(x_k)\leq \phi^\infty(\bar x)
\end{equation} for all $\bar x\in \dom \phi$. Now let $\bar x\notin \dom \phi$,      take $\hat x\in \dom \phi$ and define $x_k:=\lambda_k \hat x+(1-\lambda_k)\bar x\to \bar x$.  Then
\begin{equation*}
\phi^\infty(\bar x)  =  \sup_{t>0} \frac{\phi(\hat x+t\bar x)-\phi(\hat x)}{t}
 \geq   \frac{\phi\left(\hat x+\left(\frac{1}{\lambda_k}-1\right)\bar x\right)-\phi(\hat x)}{\frac{1}{\lambda_k}-1}
=  \lambda_k\cdot\frac{\phi\left(\hat x+\left(\frac{1}{\lambda_k}-1\right)\bar x\right)-\phi(\hat x)}{1-\lambda_k}
\end{equation*}
for all $k\in\bN$ sufficiently large. Hence for such $k\in \bN$,
\[
  (\lambda_k\star \phi)(x_k)=\lambda_k \phi\left(\frac{\lambda_k \hat x+(1-\lambda_k)\bar x}{\lambda_k}\right)\leq (1-\lambda_k)\phi^\infty(\bar x)+\lambda_k\phi(\hat x).
\]
Take the limit superior to obtain~\eqref{eq:limsup} here.
This establishes epi-convergence. 
\hfill\Halmos\endproof
 
The following result summarizes key properties of the perspective map. It also provides a support-function representation, which means that it can be written as the support function $\sigma_\cD(y)\equiv\delta^*_\cD(y)=\sup_{x\in\cD}\ip{x}{y}$
for some set $\cD$.


\begin{proposition}[Properties of perspective map]\label{prop:OmegaP}
For $\omega\in \Gamma_0(\bE_\omega)$, the following hold:
\begin{enumerate}
\item  $\omega^\pi(z,\lambda)=\sig_{\epi \omega^*}(z,-\lambda)$, hence  $\omega^\pi\in \Gamma_0(\bE_\omega\times\R)$ is sublinear with $\dom \omega^\pi=\R_+(\dom \omega\times \{1\})$;
\item $(\omega^\pi)^*(y,\beta)=\delta_{\epi\omega^*}(y,-\beta)$;
\item for all $(z,\lambda) \in \dom \omega^\pi$,
 \begin{equation}\label{eq:subdiff-omega-pi}
 \p \omega^\pi(z,\lambda)
 =\begin{cases}
	\set{(y,-\beta)}{y\in \p\omega(z/\lambda),\ \beta=\omega^*(y)} & \mbox{if $\lambda>0$,}
	\\[6pt]
	\set{(y,-\beta)}{y\in\p\omega^\infty(z),\ (y,\beta)\in\epi\omega^*} & \mbox{if $\lambda=0$.}
\end{cases}
\end{equation}
\end{enumerate}
\end{proposition}
\proof{Proof.}  For Parts (a) and (b) see~\cite[Corollary~13.5.1]{Roc 70}. Part (c) follows from \cite[Proposition~2.3]{Com 18} or \cite[Lemma 3.8]{ABD 18}.
 \hfill\Halmos\endproof

The expression for the subdifferential~\cref{eq:subdiff-omega-pi}, evaluated at the origin, reduces to $\p\omega^\pi(0,0)=\{(y,-\beta)\in\epi\omega^*\}$, which is just the epigraph of $\omega^*$ under the reflection $(z,\lambda)\mapsto(z,-\lambda)$. This follows because the subdifferential formula $\p\omega^\infty(0)=\p\sigma_{\dom\omega^*}(0)=\dom\omega^*$; cf.~\cite[Corollary~8.25]{RoW 98}. Combettes~\cite[Corollary~2.5]{Com 18} provides a simplified characterization of \cref{prop:OmegaP} under the additional assumption that $\omega$ is supercoercive~\cite[Definition~11.11]{BaC 17}.
%

\section{Partial infimal projection with perspective maps}\label{sec:Perspective}

Our main objective in this section is to deduce the variational properties of the generalized infimal convolution $\pLwf$ defined by \cref{eq:PerOpt}. Throughout this section, we make the assumptions that $L$ is a linear map from $\bE_f\times\bE_x$ to $\bE_\omega$ for Euclidean spaces $\bE_i$, $i\in\{f,x,\omega\}$, that $f\in\Gamma_0(\bE_f)$ and $\omega\in\Gamma_0(\bE_\omega)$, and that $\range L\subseteq\R_+\dom\omega$. Under these standing assumptions, it follows from \cref{th:InfProj}  below that $\pLwf$ is convex.

\subsection{Infimal projection}\label{sec:Prelim} We lead with a general result on infimal projections.


\begin{theorem}[Conjugate and subdifferentials of infimal projection]\label{th:InfProj}
For a function $\psi\in \Gamma_0(\bE_1\times \bE_2)$, the infimal projection
\begin{equation}\label{eq:inf-proj}
	p:x\in\bE_1\mapsto\inf_{u}\, \psi(x,u)
\end{equation}
is convex and 
\begin{enumerate}
\item $p^*=\psi^*(\cdot,0)$,  which is  closed and  convex;
\item $\p p(x)=\set{v}{(v,0)\in \p \psi(x,\bar u)}$ for all $\bar u\in \argmin \psi(x,\cdot)$;
\item $p^*\in \Gamma_0(\bE_1)$ if and only if
$\dom \psi^*(\cdot,0)\neq\emptyset$;
\item  $p\in \Gamma_0(\bE_1)$ if $\dom \psi^*(\cdot,0)\neq\emptyset$, and hence the infimum in \cref{eq:inf-proj} is attained when finite.
\end{enumerate}
\end{theorem}
\proof{Proof.} For convexity of $p$ and Parts (a,b,d,e), see, e.g., \cite[Theorem~3.101]{Hoh 19}. Part (c) follows from Part (b) via Rockafellar~\cite[Theorem~23.5]{Roc 70}.
\hfill\Halmos\endproof


\subsection{Generalized infimal convolution}

The following auxiliary result is used in this section to derive conjugate and a subdifferential formulas for the value function $\pLwf$.

\begin{lemma}[Domain and conjugate of linear-perspective composition]\label{lem:Eta}
  The function
  \[
	\eta:(u,x,\lambda)\in\bE_f\times\bE_x\times\R\mapsto\omega^\pi(L(u,x),\lambda)
  \]
  is closed proper convex, i.e., $\eta\in\Gamma_0(\bE_f\times\bE_x\times\R)$. The nonempty domain and its (possibly empty) relative interior are given by
  \begin{align*}
	\dom\eta\phantom)
	&= \set{(u,x,\lambda)}{\lambda \geq 0,\; L(u,x)\in \lambda \cdot \dom \omega},
  \\\ri(\dom\eta)
    &=\set{(u,x,\lambda)}{\lambda > 0,\; L(u,x)\in \lambda\cdot\ri(\dom\omega)}.
  \end{align*}
  If $\ri(\dom\eta)$ is nonempty, then $\eta^*$ is the indicator to the set
  \begin{equation}\label{eq:conjugate-domain}
	C=\set{(w,z,\mu)}{\exists y\mid (y,-\mu)\in\epi\omega^*,\; L^*(y)=(w,z)}.
  \end{equation}
\end{lemma}
\proof{Proof.} 
\Cref{prop:OmegaP}(a) asserts that $\eta\in\Gamma_0(\bE_f\times\bE_x\times\R)$, and also yields the expression for its domain. Now assume that $\ri(\dom\eta)$ is nonempty, and that there exists an element $(u,x)$ such that $L(u,x)\in\lambda\cdot\ri(\dom\omega)$ for some $\lambda>0$. Define the linear map $\tilde L:(u,x,\lambda)\mapsto (L(u,x), \lambda)$. Then,
\begin{eqnarray*}
\emptyset & \neq  & \set{(u,x,t)}{ t>0, \;L(u,x)\in t\cdot \ri(\dom\omega)}\\
& = &  \set{(u,x,\lambda)}{ \exists t>0: \;L(u,x)\in t\cdot \ri(\dom\omega),\ \lambda=t}\\
& = & \tilde L^{-1}\R_{++}(\dom \omega \times \{1\})\\
& \overset{\rm(i)}= & \tilde L^{-1} \ri(\R_+(\dom \omega \times \{1\}))\\
& \overset{\rm(ii)}= &  \ri(\tilde L^{-1}\R_+(\dom \omega \times \{1\}))\\
& = & \ri (\tilde L^{-1} \dom \omega^\pi)
= \ri(\dom\eta),
\end{eqnarray*}
where (i) uses~\cite[Corollary~6.8.1]{Roc 70} and (ii) uses~\cite[Theorem~6.7]{Roc 70} and the fact that $L^{-1}\ri(\R_+(\dom \omega \times \{1\}))\neq \emptyset$.

To derive the formula for $\eta^*$, observe that  by our reasoning above $ \tilde L^{-1}\ri(\dom\omega^\pi)=\ri(\dom\eta)\neq \emptyset$. Hence, by~\cite[Theorem~16.3]{Roc 70} and \cref{prop:OmegaP}(b),
\begin{eqnarray*}
\eta^*(w,z,\mu) & = & (\omega^\pi\circ \tilde L)^*(w,z,\mu)\\
& = & \inf_{(u,\alpha)} \set{(\omega^\pi)^*(u,\alpha)}{\tilde L^*(u,\alpha)=(w,z,\mu)}\\
& = & \inf_u \set{(\omega^\pi)^*(u,\mu)}{ L^*(u)=(w,z)}\\
& = & \inf_u\set{\delta_{\epi\omega^*}(u,-\mu)}{L^*(u)=(w,z)}\\
& = & \delta_C(w,z,\mu),
\end{eqnarray*}
which establishes~\eqref{eq:conjugate-domain}
\hfill\Halmos\endproof

We can now deduce the subdifferential and conjugate of the generalized convolution~\eqref{eq:PerOpt}.

\begin{theorem}[Conjugate and subdifferential of the generalized convolution]\label{th:ConjGProx} Under the assumptions of \cref{lem:Eta}, suppose in addition that 
\begin{equation}\label{eq:CQ1}
\exists (u,x)\in \ri(\dom f)\times \bE_x:\; L(u,x)\in \R_{++}\ri(\dom\omega).
\end{equation}
Then the following hold for the convex function $\pLwf$ defined in~\eqref{eq:PerOpt}.
\begin{enumerate}
\item   
$
\pLwf^*(y,\mu)=\inf_w \set{f^*(w)}{\exists a:(a,-\mu)\in\epi\omega^*,\ L^*(a)=(-w,y)}
$
and the infimum is attained when finite.

\item For all $(x,\lambda)\in \dom \pLwf$ and all  $\bar u\in \argmin_{u\in \bE_f}\{\,f(u)+\omega^\pi(L(u,x),\lambda)\,\}$,
\vspace{0.1cm}
\begin{equation*}
 \p \pLwf(x,\lambda)
 =
 \begin{cases}
	\set{(v,-\omega^*(y))}{ y\in \p \omega\left(L(\bar u,x)/\lambda\right),\ (0,v)\in\cD(\bar u,y) } &\mbox{if $\lambda>0,$}
	\\[6pt]
    \set{(v,-\beta)}{y\in\p\omega^\infty(L(\bar u,x)),\ (0,v)\in\cD(\bar u,y),\ (y,\beta)\in\epi\omega^*} & \mbox{if $\lambda=0$,}
\end{cases}
\end{equation*}
where $\cD(u,y):=\p f(u)\times\{0\}+L^*(y)$.
\item $\pLwf^*\in \Gamma_0(\bE_x\times \R)$ if and only if there exist $w\in \dom f^*, a\in \dom \omega^*, (y,\mu)\in \bE_x\times \R $ such that 
$(a,-\mu)\in\epi\omega^*$ and $L^*(a)=(-w,y)$. In this case, $\pLwf\in \Gamma_0(\bE_x\times \R)$ and the infimum is attained when finite.
\end{enumerate}
\end{theorem}
\proof{Proof.} Set $p=\pLwf$. Part (a). Observe that $p(x,\lambda)=\inf_u\psi(u,x,\lambda)$ for $\psi=\phi+\eta$ with $\phi(u,x,\lambda)=f(u)$ (and $\eta$ as in \cref{lem:Eta}).
We hence compute
\begin{eqnarray*}
p^*(y,\mu) & = & \psi^*(0,y,\mu)\\
& = & (\phi+\eta)^*(0,y,\mu)\\
& = & \inf_{(w,z,\delta)} \phi^*(w,z,\delta)+\eta^*(-w,y-z,\mu-\delta)\\
&= & \inf_{w} f^*(w)+\delta_C(-w,y,\mu)\\
&= & \inf_w \set{f^*(w)}{\exists a:(a,-\mu)\in\epi\omega^*,\ L^*(a)=(-w,y)}.
\end{eqnarray*}
Here the first identity uses \cref{th:InfProj}.   The second is clear from our definitions above. The third relies on~\cite[Theorem~16.4]{Roc 70} and the fact that assumption~\eqref{eq:CQ1} is, in view  of \cref{lem:Eta}(b) and the fact that $\ri(\dom\phi) = \ri(\dom f)\times\bE_x\times\R$, equivalent to the condition $\ri(\dom\eta)\cap \ri(\dom\phi)\neq \emptyset$.  The fifth uses the 
 fact that  $\phi^*(v,y,\mu)=f^*(v)+\delta_{\{0\}}(y,\mu)$ and \cref{lem:Eta}  b). The last identity is simply the definition of the set $C$ in said proposition. 

Part (b).  By~\eqref{eq:CQ1} we  can apply \cite[Theorems  23.8-23.9]{Roc 70}   to find 
\begin{align*}
  \p \psi(u,x,\lambda)
   &= \p f(u)\times \{0\}\times \{0\}+\tilde L^*\p \omega^\pi(\tilde L(u,x,\lambda))
 \\&= \p f(u)\times \{0\}\times \{0\}+(L^*\times \id)\p\omega^\pi(L(u,x),\lambda).
\end{align*}
Apply  \cref{prop:OmegaP}(c) and combine with \cref{th:InfProj} to obtain the desired result.

Part(c) follows from \cref{th:InfProj}(d).
 \hfill\Halmos\endproof

\subsection{Infimal convolution} \label{sec:InfConv}

We now consider the value function
\begin{tcolorbox}
\begin{equation}\label{eq:InfConv}
 \pwf:(x,\lambda)\in\bE\times\bR\mapsto\inf_{u\in\bE} f(u)+\omega^\pi(x-u,\lambda),
\end{equation}
\end{tcolorbox}
which corresponds to the standard infimal convolution between $f$ and $\omega^\pi$.
This is a special case of~\eqref{eq:PerOpt} where $L(u,x)= x-u$  and $\bE_i=\bE$ with $i=f,x,w$. The following result specializes \cref{th:InfProj}.

\begin{corollary}[Conjugate and subdifferential of infimal convolution]\label{cor:InfConv}
	For the function $\pwf$ given by~\eqref{eq:InfConv}, assume that $f,\omega\in\Gamma_0(\bE)$ and
\begin{equation}\label{eq:ICCQ}
\exists (u,x)\in \ri(\dom f)\times  \bE: \;x-u\in \R_{++}\ri(\dom\omega).
\end{equation}
Then the following hold.
\begin{enumerate}
\item $\pwf^*(y,\mu)=f^*(y)+\delta_{\epi \omega^*}(y,-\mu)$.
\item For all $(x,\lambda)\in \dom \pwf$ and all $\bar u\in \argmin_{u\in \bE}\left\{ f(u)+\omega^\pi(x-u,\lambda)\right\}$ we have 
\begin{equation*}
\p \pwf(x,\lambda)=
\begin{cases}
\set{(y,-\beta)}{y\in \p f(\bar u)\cap \p\omega\left(\frac{x-\bar u}{\lambda}\right),\ \beta=\omega^*(y)} & \mbox{if $\lambda>0$,}\\[6pt]
\set{(y,-\beta)}{y\in \p f(\bar u)\cap \p w^\infty(x-\bar u),\ (y,\beta)\in \epi \omega^*}&\mbox{if $\lambda =0$.}
\end{cases}
\end{equation*}

\item  $\pwf^*\in \Gamma_0(\bE)$ if and only if $\dom \pwf^*=(\dom f^*\times \bE)\cap \epi \omega^*\neq \emptyset$. In this case, $\pwf\in \Gamma_0(\bE)$ also, and the infimum is attained when finite.
\end{enumerate}

\end{corollary}
\proof{Proof.} Use \cref{th:ConjGProx}(a)--(c) and observe that $L^*(a)=(-a,a)$.
 \hfill\Halmos\endproof


\subsubsection{Infimal convolution solution map}\label{sec:Solution}

Thus far, our analysis has focused exclusively on the variational properties of the optimal value function \cref{eq:PerOpt} and its specializations. We now turn our attention to the  optimal solution map
\begin{tcolorbox}
\begin{equation}\label{eq:SolMap}
\Pwf:(x,\lambda)\in\bE_x\times\bR\mapsto\argmin_{u\in\bE_f} f(u)+\omega^\pi(x-u,\lambda)
\end{equation}
\end{tcolorbox}
for the infimal convolution defined by~\eqref{eq:InfConv}. In this section we describe the variational-analytic properties of the solution map, including (Lipschitz) continuity and (directional) smoothness. To this end, we introduce required technical machinery from variational analysis~\cite{RoW 98, Mor 18}.

Let $S:\bE_1\rightrightarrows\bE_2$ be a set-valued map between spaces $\bE_1$ and $\bE_2$. The domain and graph of $S$, respectively, are the sets $\dom S:=\set{x}{S(x)\neq \emptyset}$ and $\gph S:=\set{(x,u)\in \bE_1\times \bE_2}{u\in S(x)}$. The {\em outer limit} of $S$ at $\bar x$  is
\[
\Limsup_{x\to \bar x} S(x):=\set{y\in\bE_x}{\exists \{x_k\}\to \bar x,  \{y_k\in S(x_k)\}\to y}.
\]
Now let $A\subset  \bE$.  The {\em tangent  cone} of $A$ at $\bar x\in A$ is 
$
T_A(\bar x):=\Limsup_{t\downarrow 0} (A-\bar x)/t.
$
The {\em regular normal cone} of $A$ at $\bar x\in A$ is the polar of the tangent cone, i.e.,
$
\hat N_A(\bar x):=\set{v}{\ip{v}{y}\leq 0\;\;\forall y\in T_{A}(\bar x)}.
$
The {\em limiting normal cone}   of $A$ at $\bar x\in A$ is 
$
N_A(\bar x):=\Limsup_{x\to \bar x} \hat N_A(x).
$
The {\em coderivative} of $S$ at $(\bar x,\bar y)\in \gph S$ is the map $D^*S(\bar x\mid\bar y):\bE_2\rightrightarrows\bE_1$ defined via
\[
v\in D^*S(\bar x\mid\bar y)(y) \IFF (v,-y)\in N_{\gph S}(\bar x,\bar y).
\]
The    {\em graphical derivative} of $S$ at $(\bar x,\bar y)$ is the map $DS(\bar x\mid\bar y):\bE_f\rightrightarrows \bE_x$ given by
\[
v\in DS(\bar x \mid \bar y)( u) \IFF (u,v) \in T_{\gph S}(\bar x,\bar y),
\] 
or, equivalently, $DS(\bar x \mid \bar y)( u)=DS(\bar x | \bar y)( u)=\Limsup_{\overset{t\downarrow 0,}{u'\to u}} \frac{S(\bar x+tu')-\bar y}{t}$ \cite[Eq.~8(14)]{RoW 98}.
The {\em strict graphical derivative} of $S$ at $(\bar x,\bar y)$ is  $D_*S(\bar x\mid\bar y):\bE_f\rightrightarrows \bE_x$ given by 
\begin{eqnarray*}
D_*S(\bar x\mid\bar y)(w)=
 \set{z}{\exists \left\{\begin{array}{l}\vspace{0.1cm}\{t_k\}\downarrow 0,\, \{w_k\}\to w,\, \{z_k\}\to z,\\
  \{(x_k,y_k)\in\gph S\}\to (\bar x,\bar y) \end{array}\right\}: z_k\in\frac{S(x_k+t_kw_k)-y_k}{t_k}}.
\end{eqnarray*}
We adopt the convention to set $D^*S(\bar x):=D^*S(\bar x\mid\bar u)$ if $S(\bar x)$ is a singleton, and  proceed  analogously for the graphical derivatives.


The above generalized derivatives possess the following definiteness properties when applied to a maximally monotone operator $T:\bE\rightrightarrows\bE$, which (by definition) satisfies the inequality
\[
  \ip{v-w}{x-y}\ge0
  \quad
  \forall (v,w)\in T(x)\times T(y),	
\]
and there is no enlargement of $\gph T$ without destroying this inequality.
Our conclusion relies on {\em Minty parameterization}. 
\begin{lemma}\label{lem:PDC}
Let $T:\bE\rightrightarrows\bE$ be maximally monotone and let $(\bar y,\bar u) \in \gph T$. Then the pair $(w,z)\in\bE\times\bE$ satisfies $\ip w z \ge 0$ if one of the following conditions hold:
 \begin{enumerate}
 \item $w\in D^*T(\bar y\mid\bar u)(z)$;
 \item $z\in D_*T(\bar y\mid\bar u)(w)$;
 \item $z\in DT(\bar y\mid\bar u)(w)$.
 \end{enumerate} 
\end{lemma}
\proof{Proof.} Part (a). See~\cite[Theorem~5.6]{Mor 18}. 

Part (b). For $z\in D_*T(\bar y\mid\bar u)(w)$ there exist  $\{z_k\}\to z,  \{t_k\downarrow 0\}, \{(y_k,u_k)\in\gph T\}\to (\bar y,\bar u)$, and $\{w_k\}\to w$ such that
\begin{equation}\label{eq:Strict}
  t_k z_k\in T(y_k+t_k w_k)-u_k\quad \forall k\in \bN.
\end{equation}
Now let $\lambda>0$ and set $J_{\lambda T}:=(\lambda T+\id)^{-1}$. By Minty parameterization~\cite[Remark 23.23]{BaC 17}, there exists $\{x_k\}$ such that
$
(y_k,u_k) = \left(J_{\lambda T} (x_k),\ (x_k-J_\lambda (x_k)/\lambda\right)
$
for all $k\in \bN$.
Combining this with~\eqref{eq:Strict} yields $ x_k+t_k(\lambda z_k+w_k)\in (\lambda T+\id)(y_k+t_kw_k)$.
Thus, as $y_k=J_{\lambda T}(x_k)$, we have
$
t_kw_k=J_{\lambda T}(x_k+t_k(\lambda z_k+w_k)-J_{\lambda T}(x_k)\;(k\in\bN).
$
Because $J_{\lambda T}$ is firmly nonexpansive~\cite[Proposition~23.8]{BaC 17} and hence 1-Lipschitz, it follows that
$
\|w_k\|\leq \|\lambda z_k+w_k\| \;(k\in\bN)$,  hence
$
\|w\|\leq \|\lambda z+w\|.
$
We infer that
$
-(\lambda/2)\|z\|^2\leq \ip{z}{w}.
$
Since $\lambda>0$ was arbitrary, letting $\lambda\downarrow 0$ gives the desired inequality.

Part (c). Follows from Part (b)  and the fact that 
$DS(\bar x\mid\bar u)(w)\subset D_*S(\bar x\mid\bar u)(w)$
for all $ w\in\bE_f$.

\hfill\Halmos\endproof
 
We record another auxiliary result. Here  we call $S:\bE_f\rightrightarrows \bE_x$ {\em proto-differentiable}  at $(\bar x,\bar u)\in \gph S$ if for any $\bar z\in DS(\bar x\mid\bar u)(\bar w)$ and any $\{t_k\}\downarrow 0$ there exist $\{w_k\}\to \bar w$ and $\{z_k\}\to\bar z$ such that 
$
z_k\in (S(\bar x+ t_kw_k)-\bar u)/t_k
$ for all $k\in \bN$.

\begin{lemma}\label{lem:ProtoDiff} Let $S:\bE_1\rightrightarrows \bE_2$ be given by $S=F+T$, where $F$ is smooth and $T$ is proto-differentiable at $(\bar x,\bar u-F(\bar x))$. Then $S$ is proto-differentiable at $(\bar x,\bar u)$.  
\end{lemma}
\proof{Proof.} Let $z\in DS(\bar x\mid\bar u)(w)$ and $\{t_k\}\downarrow 0$. Then   
$
z-F'(\bar x)w\in DT(\bar x\mid\bar u-F(\bar x))(w),
$
cf.   \cite[Exercise 10.43]{RoW 98}.
By assumption on $T$, there exist $\tilde z_k\to z-F'(\bar x)w$ and $w_k\to w$ such that 
$
\tilde z_k\in [T(\bar x+t_kw_k)-(\bar u-F(\bar x))]/t_k,
$
i.e.,
$
\tilde z_k+[F(\bar x+t_kw_k)-F(\bar x)]/t_k\in[S(\bar x+t_kw_k)-\bar u]/t_k
$
for all $k\in\bN$. 
Therefore,
$
z_k:=\tilde z_k+[F(\bar x+t_kw_k)-F(\bar x)]/t_k\to z$ and  $z_k\in[S(\bar x+t_kw_k)-\bar u]/t_k$ for all $k\in \bN$
which  shows the proto-differentiability of $S$ at $(\bar x,\bar u)$.
 \hfill\Halmos\endproof

The next and main result in this subsection is based on  the implicit mapping framework described by Rockafellar and Wets~\cite[Theorem~9.56]{RoW 98} together with \cref{lem:PDC}. 

\begin{theorem}[Variational properties of the solution map]\label{th:InfConvSol}
	
	Let $f\in \Gamma_0(\bE)$ and let $\omega:\bE\to \R$  be strictly convex, level-bounded and  twice continuously differentiable. Let $\bar x\in\bE$ and $\bar \lambda>0$,  set $\bar y:=\Pwf(\bar x,\bar\lambda)$ and  $\bar V:=\nabla^2\omega\left(\frac{\bar x-\bar y}{\bar \lambda}\right)$. Then for the solution map 
$\Pwf$ from~\eqref{eq:SolMap} the following hold: 
\begin{enumerate}
\item We have $\dom \Pwf\subset\bE\times \R_{+}$ and  $\Pwf(\cdot,\lambda)$ is   single-valued for all $\lambda>0$.
\item   If $\bar V$ is positive definite, then $\Pwf$ is locally Lipschitz at $(\bar x,\bar \lambda)$.
\item If $\bar V$ is positive definite and 
$\p f$ is  {\em proto-differentiable} at $\left(\bar y,\nabla \omega\left(\frac{\bar x-\bar y}{\bar \lambda}\right)\right)$, then $\Pwf$ is  is directionally differentiable\footnote{In fact,  {\em semidifferentiable}  at $(\bar x,\bar \lambda)$ in the sense of \cite[p.~332]{RoW 98}. } at $(\bar x,\bar \lambda)$. Concretely, for all $(d,\Delta)\in \bE\times \R$, we have
\begin{eqnarray*}
 \Pwf'((\bar x,\bar \lambda);(d,\Delta))  = 
 \left[\bar \lambda D(\p f)\left(\bar y\ \Big\vert\ \nabla \omega\left(\frac{\bar x-\bar y}{\bar \lambda}\right)\right)+\bar V\right]^{-1}\hspace{-0.1cm}\left(\bar Vd-\frac{\Delta}{\bar \lambda}\bar V(\bar x-\bar y)\right).
\end{eqnarray*}
\end{enumerate}
\end{theorem}

\proof{Proof.} Set $P:=\Pwf$. Part (a). For $\lambda>0$ and $x\in\bE$, the function $u\mapsto f(y)+\omega^\pi(x-y,\lambda)$ is lsc, proper, strictly convex and level-bounded, and therefore attains a unique minimum.  

Part (b). Without loss of generality, let $\bE=\R^n$, and  observe that, for $\lambda>0$, we have  
$
P(x,\lambda)=\set{y}{0\in S(x,\lambda,y)},
$
where  $S(x,\lambda,u):= \p f(y)-\nabla \omega\left(\frac{x-u}{\lambda}\right)\;(\lambda>0)$.   Use~\cite[Exercise 10.43]{RoW 98} to deduce 
\begin{eqnarray*}
D^*S(\bar x,\bar \lambda,\bar y \mid 0)(y) =
\left[-\frac{1}{\bar \lambda} \bar Vy, \frac{(\bar x-\bar y)^T}{\bar \lambda^2}\bar Vy,\frac{1}{\bar \lambda} \bar Vy\right]+\{0\}\times \{0\}\times D^*(\p f)\left(\bar y\ \Big\vert\ \nabla \omega\left(\frac{\bar x-\bar y}{\bar \lambda}\right)\right)(y).
\end{eqnarray*}
Hence, $(r,\gamma,0)\in D^*S(\bar x,\bar \lambda,\bar y\mid 0)(y)$ if and only if
\begin{eqnarray*}
 r=-\frac{1}{\bar \lambda}\bar Vy,
 \quad \gamma=\frac{(\bar x-\bar y)^T}{\bar \lambda^2}\bar Vy,
 \quad -\frac{1}{\bar \lambda}\bar Vy\in D^*(\p f)\left(\bar y\Big\vert\nabla \omega\left(\frac{\bar x-\bar y}{\bar \lambda}\right)\right)(y).
\end{eqnarray*}
Invoke \cref{lem:PDC}(a) and use $\bar V\succ 0$ to deduce $y=0$, hence $r=0$ and $\gamma=0$. 
Therefore, by \cite[Theorem~9.56 (a)]{RoW 98}, we see that $P$ has the {\em Aubin property} at $(\bar x,\bar \lambda)$ for $\bar y$, and since $P$ is single-valued, it is locally Lipschitz at $(\bar x,\bar \lambda)$.

%

Part (c).  With the definitions from Part (b), recall that the implication 
\[
(r,\gamma,0)\in D^*S(\bar x,\bar \lambda,\bar y \mid 0)(y) \quad \Rightarrow \quad (r,\gamma)=0,\;y=0
\] 
was proved.  Now let 
\[
0\in D_*S(\bar x,\bar \lambda,\bar y\mid 0)\left(\begin{smallmatrix}0\\0\\w
\end{smallmatrix}\right)=\frac{1}{\bar \lambda}\bar Vw+D_*(\p f)\left(\bar y\mid \nabla \omega\left(\frac{\bar x-\bar y}{\bar \lambda}\right)\right)(w),
\]
see~\cite[Exercise 10.43]{RoW 98}, i.e.,
\[
-\frac{1}{\bar \lambda}\bar Vw\in  D_*(\p f)\left(\bar y\mid \nabla \omega\left(\frac{\bar x-\bar y}{\bar \lambda}\right)\right)(w).
\]
By \cref{lem:PDC}(b), we find that $w=0$. Since $\p f$ is assumed to be proto-differentiable  at $\left(\bar y,\nabla \omega\left(\frac{\bar x-\bar y}{\bar \lambda}\right)\right)$,  \cref{lem:ProtoDiff} yields that $S$ is proto-differentiable at $((\bar x,\bar \lambda,\bar y), 0)$. We can now apply \cite[Theorem~9.56(c)]{RoW 98} to obtain the desired result.
%
 \hfill\Halmos\endproof

\begin{remark}[Proto-differentiability of $\p f$ from full amenability]\label{rem:proto-diff} Let  $f\in \Gamma_0(\bE)$ and $\bar x\in \dom f$.  By \cite[Corollary~13.41]{RoW 98}, there exists a neighborhood $V$ of $\bar x$ such that   $\p f$ is proto-differentiable at $x\in V\cap \dom f$ for any $ v\in \p f(x)$ 
 if $f$ is {\em fully amenable}  at $\bar x$ in the sense that (on a neighborhood of $\bar x$)  $f=g\circ F$ with $g\in \Gamma_0(\bE_x)$ piecewise linear-quadratic and $F\in C^2(\bE_f,\bE_x)$  such that
\[
\ker F'(\bar x)^*\cap N_{\cl(\dom g)}(F(\bar x))=\{0\}.
\] 
This comprises the following special cases:
\begin{itemize}
\item $f(x)=\max_{i=1}^m f_i(x)$ with $f_i\in \Gamma_0(\bE)\cap C^2$;
\item $f$ is (convex and) piecewise linear quadratic;
\item $f$ is (convex and) twice continuously differentiable.
\end{itemize}
\end{remark}
Since  a strongly convex function is both strictly convex and level-bounded (in fact supercoercive) and has positive definite Hessian everywhere, and since we have $D(\p f)=\nabla^2 f$ wherever $f$ is twice continuously differentiable, we immediately obtain the following result  which, of course, can also be derived directly from the implicit function theorem.

\begin{corollary}[Differentiability of the solution map]\label{cor:C2case} Let $(\bar x,\bar \lambda)\in \bE\times \R_{++}$ such that  $f\in\Gamma_0(\bE)$ is  twice continuously differentiable around  $\Pwf(\bar x,\bar \lambda)$, and let $\omega\in\Gamma_0(\bE)$ be strongly convex and twice continuously differentiable. Then  $\Pwf$ from~\eqref{eq:SolMap} is continuously differentiable around $(\bar x,\bar \lambda)$. Concretely, for   all $(x,\lambda)$ sufficiently close to $(\bar x,\bar \lambda)$  and for all  $(d,\Delta)\in\bE\times \R$,  we have  
\begin{eqnarray*}
\Pwf'(x,\lambda)(d,\Delta)=\left(\lambda\nabla^2f(y)+V\right)^{-1}\left[Vd -\Delta\cdot V\left(\frac{x-y}{\lambda}\right)\right],
\end{eqnarray*}
where   $ y:=\Pwf(x,\lambda) $ and $V:=\nabla^2\omega\left(\frac{x-y}{\lambda}\right)$.
\end{corollary} 

\subsubsection{Semismoothness*} 

We now refine our study of smoothness properties of the solution map $\Pwf$. We base our analysis on the notion of {\em semismoothness*} recently established by Gfrerer and Outrata \cite{GfO 19}, which, in turn, relies on the  notion of the {\em directional normal cone}   introduced by Ginchev and Mordukohovich~\cite{GiM 12} and further advanced  by  Gfrerer et al. \cite{Gfr 13.1, Gfr 13.2, BGO 19}.

For $\bar x\in A\subset \bE$, the directional normal cone in the direction $\bar u\in\bE$ is given by 
\[
N(\bar x;\bar u)=\Limsup_{u\to \bar u,\; t\downarrow 0} \hat N_{A}(\bar x+tu).
\]
Note that $N(\bar x;\bar u)=\emptyset$ if $\bar u\notin T_A(\bar x)$ and that $N(\bar x;\bar u)\subset N_A(\bar x)$ for all $u\in \bE$. Given a set-valued map $S:\bE_f\rightrightarrows \bE_x$, based on the directional normal cone, we  define the {\em directional coderivative}  \cite{Gfr 13.1} $D^*S((\bar x,\bar u); (u,v)):\bE_x\rightrightarrows \bE_f$ of $S$ at $(\bar x,\bar y)\in\gph S$ in the direction $(u,v)$ via
\[
\gph D^*S((\bar x,\bar u); (u,v))(v^*)=\set{u^*\in\bE_f}{(u^*,-v^*)\in N_{\gph S}((\bar x,\bar y); (u,v))}.
\]
As $N(\bar x;\bar u)=\emptyset$ if $\bar u\notin T_A(\bar x)$, we also have 
\begin{equation}\label{eq:DirCoD}
\dom D^*S((\bar x,\bar u); (u,v))=\emptyset\quad  \forall (u,v)\notin DS(\bar x\mid\bar u).
\end{equation}

\begin{definition}[Semismothness*]
	The set $A\subset \bE$ is {\em semismooth*} at $\bar x\subset A$ if  
	\[
	\ip{x^*}{u}=0\quad \forall u\in \bE,\;  x^*\in N_A(\bar x;u).
	\]
	The map $S:\bE_1\rightrightarrows \bE_2$ is {\em semismooth*} at $(\bar x,\bar y)\in \gph S$ if $\gph S$ is semismooth* at $(\bar x,\bar y)$, i.e.,
	\[
	\ip{u}{u^*}=\ip{v}{v^*}\quad \forall (u,v)\in \bE_1\times \bE_2,\; (v^*,u^*)\in\gph D^*S((\bar x,\bar u); (u,v)). 
	\]
\end{definition}

\noindent
The notion of  {\em metric (sub)regularity} is  used only in the next  two results, and hence we refer the reader to the abundant literature  for a definition, e.g., \cite{DoR 14}.

\begin{proposition}[Metric regularity and semismoothness*]\label{prop:Semi} Let $F:\bE_1\to \bE_2$ be continuously differentiable at $\bar x$,  let $Q\subset \bE_2$ be semismooth*  (as a set) at $F(\bar x)$ and let $S:\bE_1\rightrightarrows \bE_2,\; S(x):=F(x)-Q$ be metrically  subregular at $(\bar x,0)$. Then $ F^{-1}(Q)$ is semismooth* at $\bar x$ (as a set).
\end{proposition}
\proof{Proof.} By \cite[Theorem~3.1]{BGO 19}, for any $h\in\bE_1$,
\begin{equation}\label{eq:DirIncl}
N_{F^{-1}(Q)}(\bar x;h)\subset F'(\bar x)^*N_{Q}(F(\bar x);F'(\bar x)h),
\end{equation}
see also \cite[Remark 2.1]{BGO 19}.  Since $Q$ is semismooth* at $F(\bar x)$,
\[
\ip{v}{z}=0\quad \forall z\in \bE_2,\ v\in N_{Q}(F(\bar x);z).
\]
Therefore
\[
\ip{v}{F'(\bar x)h}=0\quad \forall h\in\bE_1,\ v\in N_{Q}(F(\bar x);F'(\bar x)h),
\]
and hence
\[
\ip{u}{h}=0\quad \forall  h\in\bE_1,\ u\in F'(\bar x)^* N_{Q}(F(\bar x);F'(\bar x)h).
\]
By~\eqref{eq:DirIncl} this implies that 
\[
\ip{u}{h}=0\quad \forall  h\in\bE_1,\ u\in N_{F^{-1}(Q)}(\bar x;h),
\]
i.e., $F^{-1}(Q)$ is semismooth* at $\bar x$.
 \hfill\Halmos\endproof

\begin{corollary}[Semismoothness* of the infimal convolution solution map]\label{cor:SemismoothP}  Let $f\in \Gamma_0(\bE)$, let $(\bar x,\bar \lambda)\in \bE\times \R_{++}$ and let $\omega$ be strongly convex and twice continuously differentiable. Then the map $\Pwf$ from~\eqref{eq:SolMap} is semismooth* at $((\bar x,\bar \lambda), \Pwf(\bar x,\bar \lambda))$ if $\p f$ is semismooth* at $\left(\Pwf(\bar x,\bar \lambda), \nabla \omega (\frac{1}{\bar\lambda}[\bar x-\Pwf(\bar x,\bar \lambda)])\right)$.
\end{corollary}
\proof{Proof.} Without loss of generality, assume $\bE=\R^n$.  Let $F:\R^n\times \R_{++}\times \R^n\to  \R^{2n}$, 
$
F(x,\lambda,z):=(z, \nabla\omega([x-z]/\lambda).
$
Then for all $x,z\in \R^n$ and $\lambda>0$, setting $V:=\nabla^2\omega([x-z]/\lambda)\succ 0$ we have
\[
F'(x,\lambda,z)=\begin{pmatrix} 0 & 0 & I\\
\frac{1}{\lambda}V & -\frac{1}{\lambda^2}V(x-z) & -\frac{1}{\lambda}V
\end{pmatrix}.
\]
Hence, $\ker F'(x,\lambda,z)^*=\{0\}$ for all $x,z\in \R^n$, $\lambda>0$. Thus, $(x,\lambda,z)\mapsto F(x,\lambda,z)-\gph \p f$ is metrically regular. As $\gph\Pwf=F^{-1}(\gph \p f)$,  if $\p f$ is semismooth* at $\left(\Pwf(\bar x,\bar x), \nabla \omega \left(\frac{\bar x-\Pwf(\bar x,\bar \lambda)}{\bar \lambda}\right)\right)=F(\bar x,\bar \lambda,\Pwf(\bar x,\bar \lambda)) $, by \cref{prop:Semi},  $\Pwf$ is semismooth* at $((\bar x,\bar \lambda), \Pwf(\bar x,\bar \lambda))$. 
 \hfill\Halmos\endproof

\Cref{cor:SemismoothP} provides a sufficient criterion for  establishing semismoothness* of the solution map $P$ on the interior of its domain. It will be a topic of future research to exploit this on a broad scale, but we can immediately state the following result for  a function $f\in\Gamma_0(\bE)$ which is either twice continuously differentiable or  {\em piecewise linear-quadratic} (PLQ) in the sense of Rockafellar and Wets~\cite[Definition 10.20]{RoW 98}.

\begin{proposition}[Semismoothness* of the subdifferential]\label{prop:PLQ}  For $f \in\Gamma_0(\bE)$, the subgradient $\p f$ is semismooth* at $(\bar x,\bar y)\in \gph \p f$ under one of the following conditions:
\begin{enumerate}
\item $f$ is twice continuously differentiable at $\bar x$;
\item $f$ is piecewise linear-quadratic (in which case $\p f$ is semismooth* on $\bE$).
 \end{enumerate}
\end{proposition}
\proof{Proof.} 

Assume condition (a) holds.  If  $f$ is twice continuously differentiable, then $D(\p f)(\bar x\mid\bar y)=\nabla^2 f(\bar x)=D^*(\p f)(\bar x\mid\bar y)$, see \cite[Example 8.43]{RoW 98}. Now  let $(u,v)\in T_{\gph \p f}(\bar x, \bar y)$, i.e., $ v\in D(\p f)(\bar x\mid\bar y)(u) =\{\nabla^2 f(\bar x)u\}$, and  let $(x^*,y^*)\in N_{\gph \p f}((\bar x,\bar y);(u,v))\subset N_{\gph \p f}(\bar x,\bar y)$, hence $ x^*\in D^*(\p f)(\bar x\mid\bar y)(-y^*)=\{-\nabla^2 f(\bar x)y^*\}$. Thus,  we have
$
\ip{(x^*,y^*)}{(u,v)}=\ip{y^*}{\nabla^2 f(\bar x)u}-\ip{\nabla^2 f(\bar x)y^*}{u}=0.
$

Now assume condition (b) holds. It follows from \cite[Proposition~12.30]{RoW 98} that $\gph \p f$ is a finite union of polyhedra. Then \cite[Proposition~3.4/3.5]{GfO 19} yields that $\gph  \p f$ is semismooth*, which gives the desired statement.
 \hfill\Halmos\endproof

\subsection{Constrained optimization}\label{sec:relax}

We now consider an application of \cref{th:InfProj} to derive the variational properties of the optimal value of the constrained optimization problem
\begin{tcolorbox}
\begin{equation}\label{eq:LevelSetOpt}
v:(x,\lambda)\in\bE_x\times\bR\mapsto\inf_{u\in\bE_f}\set{f(u)}{L(u,x)\in \lambda S},
\end{equation}
\end{tcolorbox}
where $S\subset\bE_\omega$ is a closed convex set. This function can be viewed as a special case of~\eqref{eq:PerOpt}, where $\omega=\delta_S$ for some closed convex set $S\subset \bE_{\omega}$. To see this, it is sufficient to note that
\[
\delta^\pi_S(z,t)=\begin{cases}
	\delta_{\lambda S}(z) & \mbox{if $\lambda>0$,}\\
    \delta_{S^\infty}(z) & \mbox{if $\lambda=0$,}\\
    +\infty & \mbox{otherwise,}
\end{cases}
\] 
and thus $L(u,x)\in\lambda S$ if and only if $\delta_S^\pi(L(u,x),\lambda)$ vanishes.
Let $S^\circ:=\set{v}{\ip{v}{s}\leq 1\;\forall x\in S}$ be the polar to the set $S$.

The following result is an immediate consequence of the general study in \cref{th:ConjGProx}.

\begin{corollary}[Conjugate and subdifferential of the constrained value function]\label{eq:LevSetGen} Let $v$ be given by~\eqref{eq:LevelSetOpt} with $S\subset \bE_{\omega}$  closed and convex,  and    assume that 
\begin{equation*}
\exists u\in \ri\dom f,\ x\in \bE_x\; : \; L(u,x)\in \bR_{++}(\ri S).
\end{equation*}
Then the following hold.
\begin{enumerate}
\item We have 
 \[
 v^*(y,\mu)=\inf_{w}\set{f^*(w)}{\exists a\in-\mu S^\circ\;:\; L^*(a)=(-w,y)}.
 \]
If $S$ is a cone then 
$
v^*(y,\mu)=\inf_w\set{f^*(w)+\delta_{\R_-}(\mu)}{(-w,y)\in L^*(S^\circ)}.
$
\vspace{0.2cm}
\item For any $(x,\lambda)\in \dom v$ and  $\bar u\in \argmin_{u}\set{f(u)}{L(u,x)\in\lambda S} $,
\begin{eqnarray*}
\p  v(x,\lambda)& = &
\begin{cases} 
\set{(v,-\sigma_S(y))}{y\in N_S(L(\bar u,x)/\lambda),\; (0,v)\in\cD(\bar u, y)}
 & \mbox{if $\lambda>0$,}\\
\set{(v,-\beta)}{\exists y\in N_{S^\infty}(L(\bar u,y))\cap (\beta S^\circ)\;: \;(0,v)\in\cD(\bar u, y) } & \mbox{if $\lambda=0$,}
\end{cases}
\end{eqnarray*}
where $\cD(u,y) := \p f(u)\times\{0\}+L^*(y)$.
If $S$ is bounded (hence compact), then
\begin{eqnarray*}
 \p  v(x,\lambda)& = & \begin{cases} 
\set{(v,-\sigma_S(y))}{y\in N_S\left(L(\bar u,x)/\lambda\right),\; (0,v)\in\cD(\bar u,y)}&\mbox{if $\lambda>0$,}\\
\set{(v,-\beta)}{\exists y\in\beta S^\circ\;: \;(0,v)\in\cD(\bar u,y)}&\mbox{if $\lambda=0$.}
\end{cases}
\end{eqnarray*}
\item We have  $v^*\in\Gamma_0(\bE_x\times \R)$ if and only if  there exist $y\in \bE_x$, $w\in\dom f^*$, $\beta\in\R$ such that $(-w,y)\in -\beta L^*(S^\circ)$. In this case, also $v\in\Gamma_0(\bE_x\times \R)$ and the infimum is attained when finite.
\end{enumerate} 
\end{corollary}
\proof{Proof.}
Part (a) follows from \cref{th:ConjGProx}(a) with $w^*=\sigma_S$. If $S$ is a cone then $w^*=\delta_{S^\circ}$.
Part (b) follows from \cref{th:ConjGProx}(b), observing that $\omega^\infty=\delta_{S^\infty}$ and that $S^\infty=\{0\}$ if $S$ is bounded, in which case $N_{S^\infty}=\bE_\omega$.
Part (c) follows from (a) and \cref{th:ConjGProx}(c).
 \hfill\Halmos\endproof

\subsubsection{Relaxed linear constraints}
As an immediate specialization of \cref{eq:LevSetGen} we obtain a result on the value function
\begin{equation}\label{eq-value-fn-linear-constraints}
v:(b,\lambda)\in\R^m\times\R\mapsto\inf_{x\in \R^n}\set{f(x)}{\|Ax-b\| \le \lambda},
\end{equation}
where $f\in\Gamma_0(\R^n)$, $A\in\R^{m\times n}$ is a matrix, and $\|\cdot\|$ is any norm in $\R^n$. Denote the associated dual norm by $\|\cdot\|^\circ$, and the corresponding unit-norm ball by $\bB$.

\begin{corollary}[Relaxed linear constraints value function]\label{cor:RLC}
	If there exists a pair $(x,\lam)\in \dom f\times\R_{++}$ such that $\|Ax-b\|<\lam$, then the following hold.
\begin{enumerate}

\item (conjugate) $v^*(y,\mu)=f^*(A^Ty)+\delta_{\mu\bB^\circ}(y)$, which is closed proper convex if and only if there exists $\beta$ and $\|y\|^\circ\le \beta$ such that $A^Ty\in \dom f^*$. In this case, $v$ is closed proper convex and the infimum is attained when finite.

\item (subdifferential) For any $(b,\lambda)\in \dom v$ and $\bar x$ that achieves the infimum in \cref{eq-value-fn-linear-constraints} (and hence $\|A\bar x-b\|\le\lambda$),
\begin{eqnarray*}  \p  v(b,\lambda)& = & \begin{cases} 
\set{(y,-\|y\|^\circ)}{y\in N_\bB\left([A\bar x-b]/\lambda\right),\; -A^Ty\in \p f(\bar x)}&\mbox{if $\lambda>0$,}\\
\set{(y,-\beta)}{\|y\|^\circ\leq\beta, \;-A^Ty\in \p f(\bar x)} & \mbox{if $\lambda=0$.}
\end{cases}
\end{eqnarray*}

\item (primal existence) For $\lambda>0$ and any $b\in \R^m$, if 
\begin{equation}\label{eq:Primal}
f^\infty(y)>0 \quad\forall y\in \ker A\setminus\{0\},
\end{equation}
then 
$ \argmin_{x}\left\{f(x)+\delta^\pi_\bB(Ax-b,\lambda)\right\} \neq \emptyset$. This holds, e.g., when $f$ is level-bounded or $\rank A=n$.

\end{enumerate}

\end{corollary}
\proof{Proof.} Part (a). The expression for the conjugate $v^*$ follows from \cref{eq:LevSetGen}(a) by observing that $L:(x,b)\mapsto Ax-b$ has adjoint $L^*:z\mapsto (A^Tz,-z)$ and that $\sigma_\bB=\|\cdot\|^\circ$. The remaining claims for Part~(a) follow from \cref{th:InfProj}.

Part (b) follows from \cref{eq:LevSetGen}(b) with the foregoing observations.
d) For $\lambda>0$ and $b\in \R^m$, the effective objective function in \cref{eq-value-fn-linear-constraints} is $\phi(x):=f(x)+\delta_{\lambda\bB}(Ax-b)$. With $\hat x$ such that $\|A\hat x-b\|\leq\lambda$, which exists by the hypothesis of the theorem, we have
\[
(\delta_{\lambda\bB}\circ(A(\cdot)-b))^\infty(x) =\sup_{\tau>0}\delta_{\lambda\bB}(A\hat x-b+\tau Ax)=\delta_{\ker A}(x),
\]
where  the second identity  uses the property that $\lambda\bB$ is bounded. With \cite[Exercise 3.29]{RoW 98} we hence find that $\phi^\infty=f^\infty+\delta_{\ker A}$, which shows, using \cite[Theorem~3.26]{RoW 98}, that $\phi$ is level-bounded if~\eqref{eq:Primal}  holds.
 \hfill\Halmos\endproof

\section{Moreau envelope and proximal map}\label{sec:prox}

In this section we outline existing and new results regarding the variational properties of the Moreau envelope and the proximal map of a closed proper convex function.  

\subsection{The Moreau envelope}

The \defterm{Moreau envelope} of $f\in \Gamma_0(\bE)$ is defined by 
\[
  e_\lambda f(x):=\min_{u\in \bE}
  \left\{ f(u)+(1/2\lambda)\|x-u\|^2\right\} \quad \forall x\in \bE,\ \lambda>0,
\]
which has a Lipschitz gradient given by $\nabla e_\lambda f(x)=\frac1\lambda(x-P_\lambda f(x))$.

The following result summarizes limiting properties of the Moreau envelope as $\lambda\downarrow 0$.

\begin{proposition}[Convergence of the Moreau envelope]\label{prop:VarConv} For $f\in \Gamma_0(\bE)$, the following hold as $\lambda \downarrow 0$:
\begin{enumerate}
\item $e_\lambda f\eto f$ and $e_\lambda f\pto f$ (in fact $e_\lambda f(x)\uparrow f(x)$ for all $x\in \bE$);
\item $\lambda f\eto \delta_{\cl(\dom f)}$;
\item $\lambda e_\lambda f(x)\to \frac{1}{2}d^2_{\cl(\dom f)}(\bar x)$ as $x\to \bar x$;
\item $\lambda \p f$  converges to $N_{\cl(\dom f)}$ {\em graphically} in the sense of \cite[Definition 5.32]{RoW 98};
\item for $x\in\dom \p f$ we have  $\nabla e_\lambda f(x)\to \argmin_{g\in\p f(x)} \|g\|$. 

\end{enumerate}
\end{proposition}

\proof{Proof.} Part (a). See, e.g.,~\cite[Theorem~1.25, Proposition~7.4]{RoW 98}.

Part (b).  By \cref{lem:EpiMultConv}(b) and~\cite[Theorem~13.3]{Roc 70}, 
$
\lambda\star f^* \eto (f^*)^\infty=\sigma_{\dom f}.
$
Wijsman's theorem~\cite[Theorem~11.34]{RoW 98} then yields
$
\lambda f=(\lambda\star f^*)^*\eto
\delta_{\cl(\dom f)}.
$

Part (c). By Part (b),  $\lambda f\eto \delta_{\cl (\dom f)}$. Hence, by~\cite[Theorem~7.37]{RoW 98},
\[
\lambda e_\lambda f=e_1(\lambda f)\cto e_1\delta_{\cl(\dom f)}=\tfrac{1}{2}d^2_{\cl(\dom f)}.
\]

Part (d).  Follows from Part (b) and Attouch~\cite[Theorem~12.35]{RoW 98}.

Part (e). See~\cite[Remark 3.32]{Att 84}.
%
\hfill\Halmos\endproof

Note that~\cref{prop:VarConv}(e)  implies that there exists $K>0$ such that 
\begin{equation}\label{eq:PFrac}
\forall \bar x\in \dom \p f\; \exists K>0\;\forall \lambda>0: \; \|P_\lambda f(\bar x)-\bar x\|\leq K \lambda.
\end{equation}
\Cref{prop:VarConv}(a) suggests the following extension of the Moreau envelope at $\lambda=0$:
\begin{tcolorbox}
\begin{equation*}
  p_f:(x,\lambda)\in\bE\times\R\mapsto
  \begin{cases}
	e_\lambda f(x) & \mbox{if $\lambda>0$,}\\
    f(x)           & \mbox{if $\lambda=0$,}\\
    +\infty        & \mbox{if $\lambda<0$.}
 \end{cases}
\end{equation*}
\end{tcolorbox}
This is exactly the value function $\pwf$ from~\eqref{eq:InfConv} with $\omega=\frac{1}{2}\|\cdot\|^2$. Hence, we may rely on our general study on infimal convolution from \cref{sec:InfConv} to understand the properties of this extension of the Moreau envelope.

\begin{corollary}[Conjugate and subdifferential of the Moreau envelope]\label{cor:Moreau} Let $f\in \Gamma_0(\bE)$. Then $p_f\in\Gamma_0(\bE\times\R)$  and
\begin{enumerate}
\item $p_f^*(y,\mu)=f^*(y)+\delta_{\epi \frac{1}{2}\|\cdot\|^2}(y,-\mu)$ and $p_f^*\in \Gamma_0(\bE\times\R)$;

\item for all  $(x,\lambda)\in\dom p_f$,
\[
\p p_f(x,\lambda)=\begin{cases}
\vspace{0.1cm}
 \left(\frac{1}{\lambda}[x-P_\lambda f(x)],\ -\frac{1}{2}\|\frac1\lambda[x-P_\lambda f(x)]\|^2\right) & \mbox{if $\lambda>0$},\\
\set{(v,\beta)}{-v\in \p f(x),\ \frac{1}{2}\|v\|^2\leq \beta} & \mbox{if $\lambda=0$}.
\end{cases}
\]
\end{enumerate}

\end{corollary}
\proof{Proof.} We are in the situation of \cref{cor:InfConv} with $\omega=\frac{1}{2}\|\cdot\|^2$. In particular,  the qualification condition~\eqref{eq:ICCQ} is trivially satisfied.
 \hfill\Halmos\endproof

\subsection{Properties of the proximal map}

We now turn our attention to the proximal map. It is straightforward to show $P_\lambda f( x)\to x$ as $\lambda\downarrow 0$ for any $x\in \dom f$. The following proposition, which generalizes this statement, can be derived from monotone operator theory~\cite[Theorem~12.37]{RoW 98}. The proof that we provide here instead relies on epigraphical convergence.


\begin{proposition}[Convergence of the proximal map]\label{prop:ProxConv} Let $f\in \Gamma_0(\bE)$ and $\bar x\in \bE$. Then 
$
\lim_{\overset{\lambda \downarrow 0,}{x\to \bar x}} P_\lambda f(x)=P_{\cl (\dom f)}(\bar x).
$

\end{proposition}
\proof{Proof.} 
Let $\{\lambda_k\}\downarrow 0$,  $\{x_k\}\to \bar x$, and $\phi_k(u):=\lambda_k f(u)+\half \|u-x_k\|^2$. Use \cref{prop:VarConv}(b) to deduce $\lambda_k f\eto \delta_{\cl(\dom f)}$. Then because $\half \|(\cdot)-x_k\|^2\cto \half \|(\cdot)-\bar x\|^2$, we obtain
$
\phi_k\eto \phi:=\delta_{\cl(\dom f)}+\half \|(\cdot)-\bar x\|^2;
$
see~\cite[Theorem 7.46 b)]{RoW 98}.  Now observe that 
$P_{\lambda_k} f(x_k)=\argmin \phi_k$ and $P_{\cl(\dom f)}(\bar x)=\argmin \phi$. Since all functions $\phi_k$  are convex and $\phi$ is level-bounded (in fact, strongly convex), the sequence $\{\phi_k\}$ is, by \cite[Exercise 7.32 c)]{RoW 98}, eventually level-bounded (in the sense of \cite[p.~266]{RoW 98}). Therefore,   we can apply~\cite[Theorem~7.33]{RoW 98}, with $\varepsilon_k=0\;(k\in \bN)$, to deduce 
$
P_{\lambda_k}f(x_k)\to P_{\cl(\dom f)}(\bar x).
$
 \hfill\Halmos\endproof

We record  the following auxiliary result. 

\begin{lemma}\label{lem:Key} Let $f\in\Gamma_0(\bE)$ and fix positive scalars $\lambda$ and $\mu$. Then for all $x\in\bE$,
\begin{equation}\label{eq:p-bound1}
\begin{aligned}
  \frac{1}{2\mu}\big(\|P_\mu f(x)&-x\|^2-\|P_\lambda f(x)-x\|^2+\|P_\mu f(x)-P_\lambda f(x)\|^2\big)\\
  &\leq  f(P_{\lambda}f(x))-f(P_\mu f(x))\\
  &\leq  \frac{1}{2\lambda}\left(\|P_\mu f( x)- x\|^2-\|P_\lambda f( x)- x\|^2-\|P_\mu f( x)-P_\lambda f( x)\|^2\right),
\end{aligned}
\end{equation}
and
\begin{equation}\label{eq:p-bound2}
\|P_{\lambda}f( x)-P_\mu f( x)\|^2\leq \frac{\mu-\lambda}{\lambda+\mu}\left(\|P_\mu f( x)- x\|^2-\|P_{\lambda}f( x)- x\|^2\right).
\end{equation}
\end{lemma} 
\proof{Proof.}  Set $P(\tau):=P_\tau f(\bar x)$ for all $\tau>0$.  To obtain the bounds in~\eqref{eq:p-bound1}, use~\cite[Eq.~7(34)]{RoW 98} to infer
\[
f(x)+\frac{1}{2\tau}\|x-\bar x\|^2-f(P(\tau))-\frac{1}{2\tau}\|P(\tau)-\bar x\|^2\geq \frac{1}{2\tau}\|x-P(\tau)\|^2\quad\forall \tau>0,\ \forall x\in \bE.
\]
For $\tau=\lambda$ and $x=P(\mu)$, we hence obtain
\[
f(P(\mu))+\frac{1}{2\lambda}\|P(\mu)-\bar x\|^2-f(P(\lambda))-\frac{1}{2\lambda}\|P(\lambda)-\bar x\|^2\geq \frac{1}{2\lambda}\|P(\mu)-P(\lambda)\|^2.
\]
Analogously, for $\tau=\mu$ and $x=P(\lambda)$, we find that
\[
f(P(\lambda))+\frac{1}{2\mu}\|P(\lambda)-\bar x\|^2-f(P(\mu))-\frac{1}{2\mu}\|P(\mu)-\bar x\|^2\geq \frac{1}{2\mu}\|P(\lambda)-P(\mu)\|^2.
\]
Combining the last two inequalities now yields~\eqref{eq:p-bound1}.

Next, use~\eqref{eq:p-bound1} to obtain
\begin{align*}
  \frac{1}{\mu}\big(\|P(\mu)-x\|^2-\|P(\lambda)&-x\|^2+\|P(\lambda)-P(\mu)\|^2 \big)
  \\&\leq  \frac{1}{\lambda} \left(\|P(\mu)-x\|^2-\|P(\lambda)-x\|^2-\|P(\lambda)-P(\mu)\|^2 \right), 
\end{align*}
or, equivalently 
\[
\left(\frac{1}{\lambda}+\frac{1}{\mu}\right)\|P(\lambda)-P(\mu)\|^2\leq \left(\frac{1}{\lambda}-\frac{1}{\mu}\right)\left(\|P(\mu)-x\|^2-\|P(\lambda)-x\|^2\right),
\]
which  is equivalent to the desired inequality~\eqref{eq:p-bound2}
\hfill\Halmos\endproof

\subsection{Proximal map extension}

\Cref{prop:ProxConv} suggests the following extension of the proximal map of $f\in\Gamma_0(\bE)$: 
\begin{tcolorbox}
\begin{equation*}
   P_f:\bE\times \R\rightrightarrows \bE,\quad
   P_f(x,\lambda):=
   \begin{cases}
	 P_\lambda f(x) & \mbox{if $\lambda>0$,}\\
     P_{\cl(\dom f)}(x) & \mbox{if $\lambda=0$,}\\
     \emptyset & \mbox{if $\lambda<0$.}
\end{cases}
\end{equation*}
\end{tcolorbox}

The next result clarifies continuity properties of the proximal map extension $P_f$.

\begin{corollary}[Lipschitz continuity of the proximal map]\label{cor:Lip}
	Let $f\in \Gamma_0(\bE)$. Then $P_f$ is continuous on $\dom P_f=\bE\times \R_+$ and is locally Lipschitz on $\inter(\dom P_f)$. If $\bar x\in \dom \p f$, then $P_f$ is \defterm{upper Lipschitz} (or \defterm{calm}) at $(\bar x,0)$, and the map $\R_+\ni\mu\mapsto P_f(\bar x,\mu)$ is locally Lipschitz at $0$, i.e., there exist positive scalars $\kappa$ and $\varepsilon$  such that
	\begin{subequations}\label{cor:Lip-bc}
	\begin{alignat}{2}
	\|P_f(\bar x,0)-P_f(x,\lambda)\|&\leq \kappa\|(\bar x-x,\lambda)\| \quad &\forall(x,\lambda)&\in B_\varepsilon(\bar x,0)\cap \dom P_f,\label{cor:Lip-b}
	\\[6pt]
	\|P_f(\bar x,\lambda)-P_f(\bar x,\mu)\|&\leq \kappa |\mu-\lambda| &\forall \lambda,\mu&\in [0,\varepsilon].\label{cor:Lip-c}
	\end{alignat}
\end{subequations}
\end{corollary}
\proof{Proof.} The continuity to the boundary of the domain follows from \cref{prop:ProxConv}.  The local Lipschitz continuity on $\inter (\dom P_f)$ follows from 
\cref{th:InfConvSol} with $\omega=\frac{1}{2}\|\cdot\|^2$.

Now assume that $\bar x\in\dom\p f$, which implies $P_f(\bar x,0)=\bar x\in \dom f$. Then for all $\lambda>0$,
\[
\|P_f(x,\lambda)-P_f(\bar x,0)\|\leq \|P_\lambda f(x)-P_\lambda f(\bar x)\|+\|\bar x-P_\lambda f(\bar x)\|\leq \|x-\bar x\|+K\lambda,
\]
where $K>0$ is given via~\eqref{eq:PFrac} and we use the property that $P_\lambda f$ is  1-Lipschitz~\cite{BaC 17}.  Set $\kappa:=\max\{1,K\}$ to obtain~\eqref{cor:Lip-b}.

Let $P:=P_f(\bar x, \cdot)$.  By~\eqref{cor:Lip-b}, there exist positive scalars $\kappa$ and $\varepsilon$ such that $\|P(\tau)-\bar x\|\leq \kappa\tau$ for all $\tau\in(0, \varepsilon]$. Hence for $\mu$ and $\lambda$ in $(0,\varepsilon]$,
\begin{equation*}
\begin{array}{rcl}
\|P(\mu)-P(\lambda)\|^2
 & \leq & \frac{\mu-\lambda}{\mu+\lambda}\left(\|P(\mu)-\bar x\|^2-\|P(\lambda)-\bar x\|^2\right)\\[6pt]
& = &  \frac{\mu-\lambda}{\mu+\lambda}\left(\|P(\mu)-\bar x\|-\|P(\lambda)-\bar x\|\right)\cdot\left(\|P(\mu)-\bar x\|+\|P(\lambda)-\bar x\|\right)
\\[6pt]
& \leq  &  \frac{\mu-\lambda}{\mu+\lambda}\kappa (\mu+\lambda)\left(\|P(\mu)-\bar x\|-\|P(\lambda)-\bar x\|\right)
\\[6pt]
& \leq & \kappa |\mu-\lambda |\cdot \|P(\mu)-P(\lambda)\|,
\end{array}
\end{equation*}
where the first inequality follows from~\eqref{eq:p-bound2} of \cref{lem:Key}, and the last inequality uses the reverse triangle inequality. Use~\eqref{cor:Lip-b} to obtain~\eqref{cor:Lip-c}.
 \hfill\Halmos\endproof

The following example shows that the assumption $\bar x\in\dom\p f$ required for \cref{cor:Lip-bc} is not redundant. 

\begin{example}[Upper Lipschitz continuity of proximal map]\label{ex:NonLip}
Consider the following two functions, both contained in $\Gamma_0(\R)$:
\begin{align*}
	f(x) &=\begin{cases} -\log {x} & \mbox{if $x>0$,}\\
					     +\infty   & \mbox{otherwise},
           \end{cases}
&   g(x) &=\begin{cases} -\sqrt{x} & \mbox{if $x\geq 0$,}\\
                         +\infty   & \mbox{otherwise}.
           \end{cases}
\\\intertext{The corresponding extended proximal maps are}
  P_f(x,\lambda)
  &= \begin{cases}
	  \half(x+\sqrt{x^2+4\lambda}) & \mbox{if $\lambda>0$,} \\
      \max\{x,0\}                       & \mbox{if $\lambda=0$,}
   \end{cases}
&
 P_g(0,\lambda)&=\left(\frac{\lambda}{2}\right)^{2/3} \quad\forall \lambda\ge0;
\end{align*}
cf. Beck~\cite[Lemma~6.5]{Bec 17} for the expression for $P_f$. Observe that $\dom f$ does not include the origin, and $|P_f(0,0)-P_f(0,\lambda)|=\sqrt{\lambda}$ for all $\lambda>0$, which is not upper Lipschitz at $(0,0)$. Next, observe that $\dom \p g$ does not include the origin, and $P_g$ is not upper Lipschitz at $(0,0)$.
\hfill $\diamond$
\end{example}

The next result  on directional differentiability of $P_f$ follows from \cref{th:InfConvSol}(c) with $\omega=\frac{1}{2}\|\cdot\|^2$.

\begin{corollary}[Directional differentiability of the proximal map]\label{cor:DirDiff}
Let $f\in \Gamma_0(\bE)$ and fix $(x,\lambda)\in\bE\times \R_{++}$. If $\p f$ is proto-differentiable at $\left(P_f(x,\lambda),\, \frac1\lambda[x-P_f(x, \lambda)]\right)$, then $P_f$ is directionally differentiable at $(x,\lambda)$ with 
\begin{eqnarray*}
	P'_f((x,\lambda);(d,\Delta))
	= \left[ \lambda D(\p f)\left(P_f(x,\lambda)\ \big\vert\  \tfrac1\lambda[x-P_f(x,\lambda)]\right)+I\right]^{-1}
	\left(d-\Delta\tfrac1\lambda[x-P_f(x,\lambda)]\right)
\end{eqnarray*}
for all $(d,\Delta)\in\bE\times \R$. In particular, for any $\lambda>0$,
\[
  (P_{\lambda} f)'(x;\cdot)=\left[\lambda D(\p f)
  \left(P_f(x,\lambda)\,\big\vert\,\tfrac1\lambda[x-P_f(x,\lambda)]\right)+I\right]^{-1}(\cdot)
\]
\end{corollary} 

\subsubsection{Semismoothness* of \boldmath{$P_f$}}

We now  establish semismoothness* of the extended proximal map $P_f$ on $\bE\times \R_{++}$. We lead with an auxiliary
result.

\begin{lemma}\label{lem:SemiShift} The map $S:\bE_1\rightrightarrows \bE_2$ is semismooth* at $(y, z- y)$ if and only is $S+\id$ is semismooth* at $(y,z)$.
\end{lemma}
\proof{Proof.} The map $S$ is semismooth* at $(y,z-y)$ if and only if 
\[
\begin{array}{rcl}
 & & v\in DS(y\mid z-y)(u), \; u^*\in D^*S((y,z-y);(u,v))(v^*)\Rightarrow \ip{u}{u^*}=\ip{v}{v^*}
\\\\
 & \Longleftrightarrow &\left\{\begin{array}{l}v+u\in D(S+\id)(y|z)(u),
\\[3pt]
 u^*+v^*\in D^*(S+\id)((y,z);(u,u+v))(v^*)\end{array}\right\}\Rightarrow \ip{u}{u^*+v^*}=\ip{u+v}{v^*}
\\\\
& \Longleftrightarrow & S+\id\;\text{semismooth* at}\;(y,z).
\end{array}
\]
Here the first equivalence is  the definition of semismoothness* and~\eqref{eq:DirCoD}. The second  uses the sum rule for the 
graphical derivative \cite[Exercise 10.43]{RoW 98} and the directional coderivative \cite[Corollary 5.3 (+ comment)]{BGO 19}, respectively,  when one summand is smooth (here the identity map). The last equivalence is
a variable change and the definition of semismoothness*  and~\eqref{eq:DirCoD} again. 
 \hfill\Halmos\endproof

\begin{proposition}[Semismoothness* of \bm{$P_f$}]\label{Prop:Semi*P}
	For $f\in \Gamma_0(\bE)$,
	\begin{enumerate}
		\item $P_f$ is semismooth* at $(x,\lambda)$ if $\p f$ semismooth* at $\left(P_f(x,\lambda),\, \frac1\lambda[x-P_f(x,\lambda)]\right)$;
		
		\item $P_{\lambda} f$ is semismooth*  at $x$ if and only if $\p f$ is semismooth* at $\left(P_{\lambda} f(x),\frac1\lambda[x-P_{\lambda} f(x)]\right)$.
	\end{enumerate}
	
\end{proposition}
\proof{Proof.}
Part (a) follows from \Cref{cor:SemismoothP} with $\omega=\frac{1}{2}\|\cdot\|^2$. For Part (b), observe that $P_{\lambda} f=(\lambda \p f+\id)^{-1}$ is semismooth* at $x$ if and only if $\lambda \p f+\id$ is semismooth* at $(P_{\lambda} f(x),\,x)$ \cite[p.~7]{GfO 19}. By \cref{lem:SemiShift}, this is the case if and only if $\lambda \p f$ is semismooth* at $(P_{\lambda} f(x),\,x-P_{\lambda} f(x))$ which, in turn,  holds if and only if $\p f$ is semismooth* at $\left(P_{\lambda} f(x),\,\tfrac1\lambda[x-P_{\lambda} f(x)]\right)$.
 \hfill\Halmos\endproof

Various papers study the semismoothness \`a la Qi and Sun  \cite{QiS 93}  of $P_f$ on $\bE\times \R_{++}$.  Most of these results, trace the semismoothness of the latter back to the semismoothness 
of the Euclidean projection onto $\epi f$. The work by Meng et al.  \cite{MSZ 05, MZGS 08} deserves explicit mention, and a good discussion of these results can be found in Milzarek's thesis \cite{Mil 16}. Bearing our applications in \cref{sec:EpiLevel}  in mind, this is somewhat of a circular strategy, and hence 
we opened up a different path  via our study in \cref{sec:Solution} on semismooth* properties of solution maps.  For a map that is locally Lipschitz at a point, semismoothness* differs from traditional semismoothness only in directional differentiability as the following result by Gfrerer and Outrata \cite[Corollary 3.8]{GfO 19} shows.

\begin{lemma}[Semismooth vs. semismooth*]\label{lem:Semi} Let $F:D\subset \bE_1\to\bE_2$ be locally Lipschitz at $x\in\inter D$.  Then the following are equivalent:
\begin{enumerate}
\item $F$ is semismooth at $x$; 
\item $F$ is semismooth* and directionally differentiable at $x$.
\end{enumerate}
\end{lemma}

This lemma gives the following immediate consequence about semismoothness of $P_f$. 

\begin{corollary}[Semismoothness of \bm{$P_f$}]\label{cor:SemiP}  Let $f\in \Gamma_0(\bE)$ and fix $(x,\lambda)\in \bE\times \R_{++}$. If $\p f$ is proto-differentiable and semismooth* at 
$\left(P_f(x,\lambda),\,\tfrac1\lambda[x-P_f(x, \lambda)]\right)$, then $P_f$ is semismooth at $( x,\lambda)$. This holds, in particular, if $f$ is PLQ or  twice continuously differentiable at $P_f( x, \lambda)$, in which case $P_f$ is continuously differentiable at $( x, \lambda)$.
\end{corollary}
\proof{Proof.} For the first statement combine \cref{cor:DirDiff}, \cref{Prop:Semi*P}, and \cref{lem:Semi}. For the second, invoke \cref{rem:proto-diff} and \cref{prop:PLQ}.
 \hfill\Halmos\endproof

Note that semismoothness* does not require directional differentiability of the function in question. However, semismoothness* is still sufficient to yield convergence of  Newton-type methods under suitable regularity conditions~\cite{GfO 19, KMP 20}.  In view of the above discussion, this is  important because  the Euclidean projector onto a closed convex set may not be directionally differentiable~\cite{Sha 94}, in which case the arguments and methods  based on (standard)  semismoothness  are invalidated.

\section{The proximal value}\label{sec:eta}

The projection onto the epigraph of a function $f\in\Gamma_0(\bE)$ requires a particular value of $\lambda$ so that the equation~\eqref{eq:root} holds. In this section we examine the variational properties of the value of the proximal map as a function of $\lambda$, i.e., the function
\begin{equation}\label{eq:prox-value-lambda}
  0<\lambda\mapsto f(P_\lambda{f}(\bar x)),
\end{equation}
where $\bar x\in\bE$ is fixed. Note that this map is not generally convex, as illustrated by the following counterexample.

\begin{example}[Nonconvexity of the proximal value]\label{ex:Nonconvex} Define $f=|\cdot|+\delta_{[-1,1]}\in \Gamma_0(\R)$. By Beck~\cite[Example 6.22]{Bec 17},
\[
  P_{\lambda} f(x)=\min\{\max\{|x|-\lambda,0\},\,1\}\cdot\sgn(x)
  \quad \forall x\in \R,\ \lambda>0.
\]
Hence, for $\bar x=2$, we obtain the nonconvex function
\[
  f(P_{\lambda}f(\bar x)) =
  \begin{cases}
     1         & \mbox{if $\lambda \in (0,1]$,}\\
     2-\lambda & \mbox{if $\lambda\in (1,2]$,}\\
	 0         & \mbox{if $\lambda>2$.}
  \end{cases}
\]
\hfill $\diamond$
\end{example}

The next result describes the monotonicity and continuity of the map~\eqref{eq:prox-value-lambda}.

\begin{corollary}[Monotonicity and continuity in \bm{$\lambda$}]\label{cor:Monotone}
	Let $f\in \Gamma_0(\bE)$ and fix $\bar x\in \bE$. Then
\begin{enumerate}
\item 
$
0<\lambda\mapsto f(P_\lambda{f}(\bar x))
$
is  decreasing (i.e., increasing as $\lambda\downarrow 0$);
\item 
$
0<\lambda \mapsto \|\bar x-P_{\lambda}f(\bar x)\|
$
is increasing; 
\item $\lim_{\lambda\to 0} f(P_\lambda{f}(\bar x))=f(P_{\cl(\dom f)}(\bar x))$.
\end{enumerate}
\end{corollary}
\proof{Proof.} Parts (a) and (b).  Let $0<\lambda<\mu$ and set $P(\lambda):=P_\lambda{f}(\bar x)$, $P(\mu):=P_\mu{f}(\bar x)$, and 
 $\delta:=\frac{1}{2}(\|P(\mu)-x\|^2-\|P(\lambda)-x\|^2)$. Then from~\eqref{eq:p-bound1} of \cref{lem:Key}, we obtain
\[
\frac{1}{\mu}\delta\leq f(P(\lambda))-f(P(\mu))\leq \frac{1}{\lambda} \delta.
\]
As  $\lambda<\mu$, this implies that $\delta \geq 0$, i.e.,  $\|P(\mu)-x\|^2\geq  \|P(\lambda)-x\|^2$, and hence 
$
f(P(\mu))\leq f(P(\lambda))$. 

Part (c).  Let $\{\lambda_k\} \downarrow 0$.  Then $p_k:=P_{\lambda_k}f(\bar x)\to p:=P_{\cl(\dom f)}(\bar x)$; see  \cref{prop:ProxConv}. It follows that 
\begin{eqnarray*}
f(p)& \geq &  \limsup_{k\to\infty} \big[  f(p_k)+\tfrac{1}{2\lambda_k}\left(\|\bar x-p_k\|^2-\|\bar x-p\|^2\right)\big]\\
& \geq &\limsup_{k\to\infty} f(p_k)\\
& \geq &\liminf_{k\to\infty} f(p_k)\\
& \geq & f(p).
\end{eqnarray*}
Here the first inequality uses that $f(p)+\frac{1}{2\lambda_k}\|\bar x-p\|^2\geq  f(p_k)+\frac{1}{2\lambda_k}\|\bar x-p_k\|^2$ for all $k\in \bN$, by definition of $p_k$.  The second is due to  $\|\bar x-p_k\|\geq \|\bar x-p\|$, by the definition of $p$ and since  $p_k\in \dom f$.  The last one is just lower semicontinuity of $f$.
%
%
%
 \hfill\Halmos\endproof
%
%

As we did with the Moreau envelope and proximal map, we define the extension of the map~\eqref{eq:prox-value-lambda} to include negative values of $\lambda$:
\begin{equation*}
	\eta^f_{\bar x}:\lambda\in\R \mapsto
	\begin{cases}
	  f(P_{\lambda}f(\bar x)) & \mbox{if $\lambda>0$,}\\
	  f(P_{\cl(\dom f)}(\bar x)) & \mbox{if $\lambda\leq 0$.}
	\end{cases} 
  \end{equation*}
We call this the \emph{proximal value function}. Observe that
\begin{equation}\label{eq:eta-as-remainder}
	\eta^f_{\bar x}(\lambda)
	=
	e_\lambda f(\bar x) - (1/2\lambda)\|\bar x-P_\lambda f(\bar x)\|^2
	\quad (\lambda>0).
\end{equation}
We use \cref{cor:Monotone} to derive the following result.

\begin{corollary}[Continuity properties of the proximal value]\label{cor:eta}
	Let $f\in \Gamma_0(\bE)$ and fix $\bar x\in \bE$. Then the following hold: 
\begin{enumerate}
\item  $\eta^f_{\bar x}$ is   decreasing,  continuous (possibly in an extended real-valued sense), and finite-valued if (and only if)  $ P_{\cl(\dom f)}(\bar x)=\bar x\in \dom f$.
\item $\eta^f_{\bar x}$  is  locally Lipschitz on  $\R_{++}$.
\item If $\bar x\in \dom \p f$, then the assertion in (b) holds on $\R$.
 \end{enumerate}
\end{corollary}
\proof{Proof.} Set $\eta:=\eta^f_{\bar x}$. Parts (a) and (b). The fact that $\eta$ is decreasing follows from \cref{cor:Monotone}(a). Now consider~\eqref{eq:eta-as-remainder}. 
By Corollary \ref{cor:Moreau}, the map $0<\lambda\mapsto e_\lambda f(\bar x)$ is convex and finite-valued, hence locally Lipschitz. By \cref{cor:Lip}(a), this conclusion also holds for $0<\lambda\mapsto \frac{1}{2\lambda}\|x-P_\lambda f(\bar x)\|^2$. This gives the local Lipschitz continuity of $\eta$ on $\R_{++}$. The continuity at 0 is due to \cref{cor:Monotone}(c).

Part (c).   By Parts (a) and (b), and because $\eta$ is constant (and finite by assumption) on $\R_-$, we only need to be concerned about 
the desired properties at $0$. To this end, let $\mu>\lambda$. If $\lambda<0$, then 
\[
\left| \frac{\eta(\mu)-\eta(\lambda)}{\mu-\lambda}\right|\leq \left|\frac{\eta(\mu)-\eta(0)}{\mu-0}\right|.
\]
Thus we can restrict ourselves to the case $0\leq \lambda<\mu$. 
Set $P(\tau):=P_{\tau}f(\bar x)$ for all $\tau> 0$ and $P(0):=\bar x$.  Then by \cref{cor:Lip}(c), there  exist positive scalars $\varepsilon$ and $\kappa$ such that 
\begin{equation}\label{eq:Lip1}
\|P(\mu)-P(\lambda)\|\leq \kappa(\mu-\lambda)\quad \forall 0\leq \lambda\leq \mu\leq \varepsilon.
\end{equation}
For $0<\lambda<\mu\leq \varepsilon$, we have
\begin{equation*}
\begin{array}{rcl}
\vspace{0.2cm}
|\eta(\lambda)-\eta(\mu)|
& = & \eta(\lambda)-\eta(\mu)\\
\vspace{0.2cm}
& \leq &\frac{1}{2\lambda}\left(\|P(\mu)-\bar x\|^2-\|P(\lambda)-\bar x\|^2-\|P(\mu)-P(\lambda)\|^2\right)\\
\vspace{0.2cm}
& = & \frac{1}{2\lambda}\left[\left(\|P(\mu)-\bar x\|-\|P(\lambda)-\bar x\|\right)\cdot\left(\|P(\mu)-\bar x\|+\|P(\lambda)-\bar x\|\right)-\|P(\mu)-P(\lambda)\|^2\right]\\
\vspace{0.2cm}
& \leq& \frac1{2\lambda}\|P(\mu)-P(\lambda)\|\cdot\left(\right \|P(\mu)-\bar x\|+\|P(\lambda)-\bar x\|-\|P(\mu)-P(\lambda)\|)\\
\vspace{0.2cm}
& \leq & \frac{\kappa}{2\lambda}|\mu-\lambda|\left(  \|P(\mu)-\bar x\|+\|P(\lambda)-\bar x\|-( \|P(\mu)-\bar x\|-\|P(\lambda)-\bar x\|)\right)\\
\vspace{0.2cm}
& = & \frac\kappa\lambda\|\bar x-P(\lambda)\|\cdot|\mu-\lambda|\\
& \leq & \kappa^2|\mu-\lambda|.
\end{array}
\end{equation*}
Here, the first identity follows from \cref{cor:Monotone}(a), where the first inequality uses \cref{lem:Key}(a). The rest the follows from the reverse triangle inequality and~\eqref{eq:Lip1}, recalling that $P(0)=\bar x$.
 \hfill\Halmos\endproof

\begin{remark} The requirement that $\bar x\in\p f$, made in \cref{cor:eta}, cannot be relaxed to $\bar x\in \dom f$. To see this, we again use \cref{ex:NonLip}(b), where
\[
  \eta^f_{\bar x}(\lambda)=
  \begin{cases}
	-\left(\lambda/2\right)^{\frac{1}{3}} & \mbox{if $\lambda\geq 0$,}\\
     0  & \mbox{if $\lambda<0$,}
\end{cases}
\]
which is neither locally Lipschitz nor directionally differentiable at $\lambda=0$. We also conclude from this example that the lack of calmness of the proximal map at $\lambda=0$ is not necessarily compensated by applying $f$.
\end{remark}

Under certain assumptions described by \cref{cor:Phi}, we may interpret the extended proximal value function $\eta^f_{\bar x}$ as the derivative of the convex function
\begin{equation}\label{eq:phibar}
	\bar \phi_{\bar x}^f:\lambda\in\R\mapsto
	\begin{cases}
		-\lambda e_\lambda f(\bar x) & \mbox{if $\lambda>0$,}\\
	    -\frac{1}{2}d^2_{\cl(\dom f)}(\bar x) & \mbox{if $\lambda=0$,}\\
		-\lambda f(P_{\cl(\dom f)}(\bar x))-\frac{1}{2}d^2_{\cl(\dom f)}(\bar x)
			& \mbox{if $\lambda<0$;}
	\end{cases}
	\end{equation}	
cf.~Attouch~\cite[Remark 3.32]{Att 84}.

\begin{corollary}[The function \bm{$\bar \phi_{\bar x}^f$}] \label{cor:Phi} Let $f\in \Gamma_0(\bE)$ and fix $\bar x\in \bE$. Then the following hold:
 \begin{enumerate}
 \item$\bar \phi_{\bar x}^f$ is proper, convex and continuous (possibly in an extended real-valued sense), and continuously differentiable on $\R_{++}$ with $\frac{d}{d\lambda}\bar \phi_{\bar x}^f(\lambda)=-f(P_{\lambda}f(\bar x))$ locally Lipschitz
   for  all $\lambda>0$.

\item If $\bar x\in \dom f$, then $ \bar \phi_{\bar x}^f$ is continuously differentiable on $\R$ with  derivative given by
\[
\frac{d}{d\lambda}\bar \phi_{\bar x}^f(\lambda)
= -\eta_{\bar x}^f(\lambda)
=
\begin{cases}
 -f(P_{\lambda}f(\bar x))& \mbox{if $\lambda>0$,}\\
 -f(P_{\cl(\dom f)}(\bar x)) & \mbox{if $\lambda \leq 0$.}
\end{cases}
\]
If, more strictly, $\bar x\in\dom \p f$, then this derivative is locally Lipschitz 
 on all of $\R$.
\item If $P_{\cl(\dom f)}(\bar x)\notin \dom f$, then $\dom \bar \phi_{\bar x}^f=\R_+$ and 
\[
\p \bar \phi_{\bar x}^f(\lambda)=
\begin{cases}
  -f(P_{\lambda}f(\bar x)) & \mbox{if $\lambda>0$,}\\
  \emptyset & \mbox{if $\lambda \leq 0$.}
\end{cases}
\]
\end{enumerate}

\end{corollary}
\proof{Proof.} Set  $\bar \phi:=\bar \phi_{\bar x}^f(\lambda)$.
Part (a).  It is an easy computation to see that 
\[
0<\lambda\mapsto -\bar \phi(\lambda)=\inf_{u}\left\{\lambda f(y)+\tfrac{1}{2}\|u-\bar x\|^2\right\}
\]
is concave, i.e., $0<\lambda \mapsto \bar \phi(\lambda)$ is convex. By setting $\bar \phi(0)=-\frac{1}{2}d^2_{\cl(\dom f)}(\bar x)$ and using \cref{prop:VarConv}(a), we see that $\bar \phi$ is a continuous convex function on $\R_+$, which is  linearly extended to $\R_-$.  All in all, $\bar \phi$
is convex, proper and continuous (possibly in an extended real-valued) sense.
From \cref{cor:Moreau}(b) (and the product rule) we infer, for all $\lambda>0$, that 
\begin{equation*}
\bar \phi'(\lambda)
 =  -e_\lambda f(\bar x)-\lambda\left(-\tfrac{1}{2}\left\|\tfrac1\lambda[\bar x-P_\lambda f(\bar x)]\right\|^2\right)
 =  (1/2\lambda)\|\bar x-P_\lambda f(\bar x)\|^2-e_\lambda f(\bar x)
 =  -f(P_\lambda f(\bar x)),
\end{equation*}
where the last equality follows from~\eqref{eq:eta-as-remainder}. Hence, the local Lipschitz continuity  follows from \cref{cor:eta}(b).

Part (b). Here we assume that $P_{\cl(\dom f}(\bar x)\in \dom f$.  Then by definition of $\bar \phi$, we have  $\bar \phi'(\lambda)=-f(P_{\cl(\dom f)}(\bar x))$ for all $\lambda<0$.
It remains to establish the case $\lambda =0$. To this end, use the subgradient inequality to deduce that $g\in \p\bar \phi(0)$ if and only if $\bar \phi(0)+\lambda g\leq \bar \phi(\lambda)$ for all $\lambda$ if and only if
\begin{subequations}
\begin{alignat}{2}
   g + e_\lambda f(\bar x) - \tfrac{1}{2\lambda}d^2_{\cl(\dom f)}(\bar x) &\le0 &\quad \forall \lambda&>0, \label{eq:subgrad-g-a}
   \\
   \lambda g + \lambda f(P_{\cl(\dom f)}(\bar x)) &\leq 0 & \forall \lambda&<0, \label{eq:subgrad-g-b}
\end{alignat}
\end{subequations}
hold simultaneously. (The case with $\lambda=0$ holds trivially.) From~\eqref{eq:subgrad-g-a}, we infer that
\begin{eqnarray*}
  g\le\inf_{\lambda >0}
  \frac{\half d^2_{\cl(\dom f)}(\bar x)-\lambda e_\lambda f(\bar x)}{\lambda}
  &\overset{\rm(i)}{=}& \inf_{\lambda >0}\frac{\bar \phi(\lambda)-\bar \phi(0)}{\lambda}
\\&\overset{\rm(ii)}{=}& \lim_{\lambda\downarrow 0} \frac{\bar \phi(\lambda)-\bar \phi(0)}{\lambda}
\\&\overset{\rm(iii)}{=}& \lim_{\lambda\downarrow 0}\frac{\half d^2_{\cl(\dom f)}(\bar x)-\lambda e_\lambda f(\bar x)}{\lambda} \\
&\overset{\rm(iv)}{=}& \lim_{\lambda\downarrow 0} \frac{-f(P_{\lambda} f(\bar x))}{1}\\
&\overset{\rm(v)}{=}& -f(P_{\cl(\dom f)}(\bar x)).
\end{eqnarray*}
Here, (i) is simply the definition of $\bar \phi$; (ii) holds because $\bar \phi$ is convex~\cite[Theorem~23.1]{Roc 70}; (iii) follows from the definition of $\bar \phi$; and (iv) follows from  l'H\^{o}pital's rule, which is applicable because the last limit exists by \cref{cor:Monotone}(c), which implies (v). Hence,~\eqref{eq:subgrad-g-a} is equivalent to $g\leq -f(P_{\cl(\dom f}(\bar x))$. Combined with~\eqref{eq:subgrad-g-b}, which is equivalent to $g\geq -f(P_{\cl(\dom f}(\bar x))$, establishes that $\p \bar \phi(0)=\{-f(P_{\cl(\dom f}(\bar x))\}$. Thus, $P_{\cl(\dom f}(\bar x)\in \dom f$, $\bar \phi$ is  differentiable, and hence continuously differentiable by convexity~\cite[Corollary 25.5.1]{Roc 70}. The remainder follows from \cref{cor:eta}(c).

Part (c). Here we assume that $P_{\cl(\dom f}(\bar x)\notin \dom f$.  Suppose $g\in \p \phi(0)$, i.e., analogous to some arguments in b),
\[
g\leq (1/2\lambda)d^2_{\cl(\dom f)} (\bar x)-e_\lambda f(\bar x)\quad \forall \lambda>0.
\]
On the other hand, using e.g., \cref{cor:Monotone}(b), we have
\begin{align*}
(1/2\lambda)d^2_{\cl(\dom f)} (\bar x)-e_\lambda f(\bar x)
 &= (1/2\lambda)\|\bar x-P_{\cl(\dom f)}(\bar x)\|^2- (1/2\lambda)\|\bar x-P_\lambda f(\bar x)\|^2-f(P_\lambda f(\bar x))\\
 &\leq  -f(P_\lambda f(\bar x)).
\end{align*}
Since $ -f(P_\lambda f(\bar x))\to  -\infty $ as $\lambda\downarrow 0$, 
this concludes the proof. 
 \hfill\Halmos\endproof

\subsection{Semismoothness of the proximal value function}  \label{sec:semismooth}
 
In view of the properties of the proximal value function, as outlined by \cref{cor:eta},   the  question for  {\em semismoothness} of $\eta_{\bar x}^f$ on $\R_{++}$ arises naturally.  Now consider the expression~\eqref{eq:eta-as-remainder}. The map $0<\lambda \mapsto e_{\lambda}f(\bar x)$ is continuously differentiable by \cref{cor:Moreau}(a), hence semismooth~\cite[Proposition~7.4.5]{FaP 03}. Moreover,
the map $0<\lambda\mapsto(1/2\lambda)\|\bar x-P_{\lambda}f(\bar x)\|^2$ is semismooth if $0<\lambda\mapsto P_{\lambda}f(\bar x)$ is semismooth~\cite[Proposition~7.4.4]{FaP 03}. Thus, when the latter holds,  we can conclude that $\eta^f_{\bar x}$ is semismooth. We can in addition use \cref{cor:SemiP}, which establishes conditions for  the semismoothness of the map $(x,\lambda)\in \bE\times \R_{++}\mapsto P_{\lambda}f(x)$, to obtain the following result. 
 
\begin{proposition}[Semismoothness of the proximal value function]\label{prop:SemiEta} Let $f\in\Gamma_0(\bE)$ and $\bar x\in\bE$. Then 
$\eta_{\bar x}^f$ is semismooth at $\bar \lambda>0$ if $\p f$ is proto-differentiable and semismooth* at 
$\left(P_{\bar\lambda}f(\bar x),\tfrac1\lambda[\bar x-P_{\bar\lambda}f(\bar x)]\right)$. This is the case   under either of the following conditions:
\begin{enumerate}
\item (PLQ case) $f$ is piecewise-linear quadratic.
\item ($C^2$ case) $f$ is twice continuously differentiable around $P_{\bar \lambda}f(\bar x)$. In this case, $\eta_{\bar x}^f$ is continuously differentiable.
\end{enumerate}
\end{proposition}

\section{Post-composition envelopes and proximal maps}\label{sec:post-composition}

Given functions $\psi\in \Gamma_0(\bE)$ and  $g\in \Gamma_0(\R)$, we  consider the composition 
\begin{equation*}\label{eq:comp}
  (g\circ \psi)(x):=\begin{cases} g(\psi(x)) & \mbox{if $x\in \dom \psi$,}\\
+\infty & \mbox{otherwise.} 
\end{cases}
\end{equation*}
It is well known that $g\circ \psi$ is closed proper convex if $g$ is increasing and that the intersection $\psi(\bE)\cap \dom g$ is nonempty; see, for example, Hiriart-Urruty and Lemar\'echal~\cite[Theorem~B.2.1.7]{HUL 01}, who describe this operation as {\em post-composition}. We establish variational formulas for the Moreau envelope  and proximal map of the composition $g\circ \psi$ under a regularity assumption involving the intersection of domains. These results provide us with tools to infer properties of projections onto the epigraph and level sets of a closed proper convex function, as covered in \cref{sec:EpiLevel}.


\begin{proposition}[Post-composition, Moreau envelopes, and proximal maps]\label{prop:PostComp} Let $g\in \Gamma_0(\R)$   be increasing and let $\psi\in \Gamma_0(\bE)$ such that \begin{equation}\label{eq:PCCQ}
(\ri\dom g) \cap \psi(\ri\dom \psi)\neq \emptyset.
\end{equation}
Then  the following properties hold.
\begin{enumerate}
\item $e_1(g\circ \psi)(\bar x)=-\min_{\lambda\geq 0}\left\{g^*(\lambda)+\bar \phi^{\psi}_{\bar x}(\lambda)\right\}$, where $\bar \phi^{\psi}_{\bar x}$ is given by~\eqref{eq:phibar}.
\item$P_1(g\circ \psi)(\bar x)=P_1(\bar \lambda\cdot \psi)(\bar x)$ for every $\bar \lambda\in\argmin_{\lambda\geq 0}\left\{ g^*(\lambda)+\bar \phi^{\psi}_{\bar x}(\lambda)\right\}\neq\emptyset$.
\item  If $\psi (P_{\cl(\dom \psi)}(\bar x))\notin\p g^*(0)$, then 
$
\argmin_{\lambda\geq 0}\left\{ g^*(\lambda)+\bar \phi^{\psi}_{\bar x}(\lambda)\right\}\subset \R_{++}$.
This is, in particular, the case if $P_{\cl(\dom \psi)}(\bar x)\notin \dom \psi$. 
 
\end{enumerate}

\end{proposition}

\proof{Proof.} Part (a). We find that
\begin{eqnarray*}
e_1 (g\circ \psi)(\bar x) & = & \min_{x\in \bE} \left\{\tfrac{1}{2}\|x-\bar x\|^2+(g\circ \psi)(x)\right\}\\
& = & {}{-}\left(\tfrac{1}{2}\|(\cdot)-\bar x\|^2+g\circ \psi\right)^*(0)\\
& = & \max_{y\in\bE,\lambda\geq 0} -\left\{g^*(\lambda)-\tfrac{1}{2}\|y\|^2+\ip{\bar x}{y}+(\lambda\cdot\psi)^*(-y)\right\}\\
& = & \max_{\lambda\geq 0}\big\{ -g^*(\lambda)+\max_{y\in\bE}\big[-\tfrac{1}{2}\|y\|^2-\ip{\bar x}{y}-(\lambda\cdot\psi)^*(-y)\big]\big\}\\
&= & \max_{\lambda\geq 0}\ -g^*(\lambda)-\bar \phi^{\psi}_{\bar x}(\lambda).
\end{eqnarray*}
Here, the third identity uses   \cite[Corollary 3]{BHN 20}  with $f:=\frac{1}{2}\|(\cdot)-\bar x\|^2$, $F:=\psi$, and $K=\R_+$, realizing that~\eqref{eq:PCCQ} is equivalent to qualification condition \cite[Equation (17)]{BHN 20} because $\dom g-\R_+=\dom g$, and observing that attainment is guaranteed by finiteness of the left-hand side. The last identity uses Fenchel duality~\cite[Theorem~31.1]{Roc 70} and the definition of $\bar \phi^{\psi}_{\bar x}$ in~\eqref{eq:phibar}.

Part (b). Note that by~\cite[Corollary 4]{BHN 20},
\begin{equation}\label{eq:SDPC} 
\p(g\circ \psi)(x)=\ \bigcup_{\mathclap{\lambda \in \p g(\psi(x))}}\ \p(\lambda\cdot\psi)(x)\quad \forall x\in \dom g\circ \psi,
\end{equation}
and observe that $\p g(x)\subset \R_+$ because $g$ is increasing. 
Next, observe that 
\begin{equation*}
\begin{array}{rcl}
\bar \lambda \in \argmin_{\lambda\geq 0} \{g^*(\lambda)+\bar \phi^{\psi}_{\bar x}(\lambda)\},\; \bar u=P_{1}(\bar \lambda\cdot\psi)(\bar x)
&\overset{\rm(i)}\Longleftrightarrow & 0\in \p g^*(\bar \lambda)+\p\bar \phi^\psi_{\bar x}(\bar \lambda),\; \bar u=P_{1}(\bar \lambda\cdot\psi)(\bar x)\\
&\overset{\rm(ii)}\Longleftrightarrow & \psi(\bar u)\in \p g^*(\bar \lambda), \; \bar u=P_{1}(\bar \lambda\cdot\psi)(\bar u)\\
&\overset{\rm(iii)}\Longleftrightarrow & \bar \lambda \in \p g( \psi(\bar u)),\; \bar u=P_{1}(\bar \lambda\cdot\psi)(\bar x)\\
&\overset{\rm(iv)}\Longleftrightarrow &  \bar \lambda \in \p g(\psi(\bar u)),\; 0\in  \bar u-\bar x+\p (\bar \lambda\cdot\psi)(\bar x)\\
&\overset{\rm(v)}\Longrightarrow &  \bar u=P_{1}(g\circ \psi)(\bar x).
\end{array}
\end{equation*}
Equivalence (i) is valid because $\inter(\dom g^*)\subset \R_{++}\subset \inter(\dom \bar \phi^{\psi}_{\bar x})$; see \cite[Lemma 4]{BHN 20}  and  \cref{cor:Phi}, respectively.  \Cref{cor:Phi}(b) justifies equivalence (ii). Equivalence (iii) is the inversion formula for the subdifferential \cite[Corollary 23.5.1]{Roc 70}. Equivalence (iv) uses the optimality conditions that uniquely determines $\bar u=P_{1}(\bar \lambda\cdot\psi)(\bar x)$. Implication (v) follows from~\eqref{eq:SDPC} and the optimality conditions that uniquely determine $P_1(g\circ \psi)(\bar x)$. Taken together, we deduce that for any $\bar \lambda\in  \argmin_{\lambda\geq 0} g^*(\lambda)+\bar \phi^{\psi}_{\bar x}(\lambda)$, we have $P_1(g\circ \psi)(\bar x)=P_1(\lambda\cdot \psi)(\bar x)$. The fact that $\argmin_{\lambda\geq 0}\left\{ g^*(\lambda)+\bar \phi^{\psi}_{\bar x}(\lambda)\right\}\neq \emptyset$ follows from Part (a).

Part (c). Recall from Part (b) that $0\in \argmin_{\lambda \geq 0}\{g^*+\phi_{\bar x}^\psi\}$ entails 
$
0\in \p g^*(0)+\p\phi_{\bar x}^f(0).
$
In view of \cref{cor:Phi}(c), we must have $P_{\cl(\dom \psi)}(\bar x)\in \dom \psi$, in which case $\p\phi_{\bar x}^\psi(0)=-\psi(P_{\cl(\dom \psi)}(\bar x))$, by \cref{cor:Phi}(b). This proves the claim.
 \hfill\Halmos\endproof

\section{Epigraphical and level-set projections}\label{sec:EpiLevel}

We are now equipped to answer the initial question about computing epigraphical and level-set projections via proximal mappings. Our approach is based on the Moreau envelopes of the indicator functions to the epigraph and level set of a function $f$, which we express as the post-compositions
\[
  \delta_{\lev{\alpha}{f}} = (\delta_{\R_-})\circ(f(\cdot)-\alpha)
  \quad\UND\quad
  \delta_{\epi f} = (\delta_{\R_-})\circ(f(\cdot)-(\cdot)).
\]
\Cref{prop:PostComp} provides the required tools.

\begin{corollary}[Level-set projection]\label{cor:LevelSet}
	Let $f\in \Gamma_0(\bE)$, $(\bar x,\,\bar \alpha)\in\bE\times\R$, and assume there exists $\hat x\in\bE$ such that $f(\hat x)<\bar \alpha$. Then the following statements hold.
	\begin{enumerate}
		\item(Dual representation of distance to level set)  
		\[
		\tfrac{1}{2}d^2_{\lev{\bar \alpha}{f}}(\bar x)=-\min_{\lambda\geq 0}\left\{ \bar \phi_{\bar x}^f(\lambda)+\bar \alpha\lambda\right\}.
		\]
		\item  (Projection onto level set)
		\[
		P_{\lev{\bar \alpha}{f}}(\bar x)=
		\begin{cases}
			P_{\cl(\dom f)}(\bar x) & \mbox{if $f(P_{\cl(\dom f)}(\bar x))\leq \bar\alpha$,}  \\
			P_{\bar \lambda}f(\bar x) & \mbox{otherwise,}
		\end{cases}
		\]
		for any positive $\bar\lambda$ in the optimal solution set
		\begin{equation*}
			\argmin_{\lambda\geq 0}\ \{\bar \phi_{\bar x}^f(\lambda)+\bar \alpha\lambda\}
			=\set{\lambda\geq 0}{f(P_{\lambda}f(\bar x))=\bar \alpha}\neq\emptyset.
		\end{equation*}
	\end{enumerate}
\end{corollary}  

\proof{Proof.}
Set $g:=\delta_{\R_-}$ and $\psi:x\in\bE\mapsto f(x)-\bar\alpha$. Then $g\in\Gamma_0(\R)$ is increasing and $\psi\in \Gamma_0(\bE)$ with $\dom \psi =\dom f$ and 
$\delta_{\lev{\bar \alpha}{f}}=g\circ \psi$. Now observe that~\eqref{eq:PCCQ} applied to this setting is equivalent to saying that there exists $\bar y\in \ri(\dom f)$ such that $f(\bar y)<\bar \alpha$. We (only) assume that there exists $\hat x\in\dom f$ such that $f(\hat x)<\bar \alpha$. However, take any $z\in \ri(\dom f)$, then, by the line segment principle~\cite[Theorem~6.1]{Roc 70},  we have $y_\lambda:=\lambda z+(1-\lambda)\hat x\in \ri(\dom f)$ for all $\lambda \in (0,1]$. Moreover,
$f(y_\lambda)< \lambda f(z)+(1-\lambda)\bar\alpha\to \bar\alpha$ as $\lambda \downarrow 0$. Hence there exists $\hat \lambda\in(0,1]$ sufficiently small such that $f(y_{\hat \lambda})<\bar\alpha$. Hence $\hat y:=y_{\hat \lambda}\in\ri(\dom f)$ with $f(\hat y)<\bar \alpha$, and~\eqref{eq:PCCQ} holds.

Part (a). For all $\lambda \geq 0$,
\begin{equation*}
\begin{array}{rcl}
\vspace{0.2cm}
\bar \phi_{\bar x}^\psi(\lambda)& = &
\begin{cases}
  -\lambda e_{\lambda}\psi(\bar x) & \mbox{if $\lambda>0$,}\\
  -\frac{1}{2}d^2_{\cl(\dom \psi)}(\bar x) & \mbox{if $\lambda=0$,}
\end{cases}\\
\vspace{0.2cm}
& =&
\begin{cases}
-\lambda( e_{\lambda}f(\bar x)-\bar \alpha) & \mbox{if $\lambda>0$,}\\
-\frac{1}{2}d^2_{\cl(\dom f)}(\bar x) & \mbox{if $\lambda=0$,}
\end{cases}\\
&= & \bar \phi_{\bar x}^f(\lambda)+\bar \alpha\lambda.
\end{array}
\end{equation*}
Use \cref{prop:PostComp}(a) and the fact that  $ g^*=\delta_{\R_+}$ to deduce that
\[
 \tfrac{1}{2}d^2_{\lev{\bar \alpha}{f}}(\bar x)
 = e_1\delta_{\lev{\bar \alpha}{f}}(\bar x)
 = e_1(g\circ \psi)(\bar x)
 = -\min_{\lambda\geq 0}\left\{ \bar \phi_{\bar x}^f(\lambda)+\bar \alpha\lambda\right\}.
\]

Part (b). The equality of the two sets in question is clear from the (necessary and sufficient) optimality conditions and \cref{cor:Phi}.   The rest follows from \cref{prop:PostComp}, Parts (b) and (c) because $P_{\lev{\bar \alpha}{f}}(\bar x)=P_1(g\circ \psi)(\bar x)$.
 \hfill\Halmos\endproof

\begin{corollary}[Epigraphical projection]\label{cor:Epi}
	Let $f\in \Gamma_0(\bE)$ and $(\bar x,\,\bar \alpha)\in\bE\times\R$. Then the following statements hold.
	\begin{enumerate}
		\item (Dual representation of distance to epigraph)
		\[
		\tfrac{1}{2}d^2_{\epi f}(\bar x,\bar \alpha)
		=-\min_{\lambda\geq 0}\left\{
		\bar \phi_{\bar x}^f(\lambda)+\bar \alpha\lambda+\tfrac12\lambda^2
		\right\}.
		\]
		\item  (Projection onto epigraph)
		\[
		P_{\epi f}(\bar x)=
		\begin{cases}
			[P_{\cl(\dom f)}(\bar x),\ \bar \alpha] & \mbox{if $f(P_{\cl(\dom f)}(\bar x))\leq \bar\alpha$,}
			\\ [P_{\bar \lambda}f(\bar x),\ \bar \alpha+\bar \lambda] & \mbox{otherwise,}
		\end{cases}
		\]
		where $\bar \lambda>0$ is the unique solution of the strongly convex optimization problem
		\begin{equation*}
			\min_{\lambda \geq 0}\ \tfrac{1}{2}\lambda^2+\bar \alpha \lambda+ \bar \phi_{\bar x}^f(\lambda).
		\end{equation*}
		Equivalently, $\lambda$ is the  unique root of the strictly decreasing function  $0<\lambda\mapsto f(P_{\lambda}f(\bar x))-\lambda-\bar \alpha$.
		
	\end{enumerate}
\end{corollary}  

\proof{Proof.} Analogous to the proof of \cref{cor:LevelSet}, we define closed proper convex functions $g:=\delta_{\R_-}$ and $\psi:(x,\alpha)\in \bE\times \R\mapsto f(x)-\alpha $ so that $\delta_{\epi f}=g\circ \psi$. Therefore,
\[
\psi(\ri(\dom \psi))=\psi(\ri(\dom f)\times \R)=f(\ri(\dom f))-\R=\R,
\]
and thus the qualification condition~\eqref{eq:PCCQ} is trivially satisfied in this setting.

Part (a). Note that $e_\lambda \psi( x,\alpha)=e_\lambda f(x)+e_\lambda (-\id)(\alpha)$ for all $\lambda>0$~\cite[Theorem~6.58]{Bec 17}, and since $\dom \psi=\dom f$,
\[
\bar \phi^\psi_{\bar x,\bar \alpha}(\lambda)=\bar  \phi^f_{\bar x}(\lambda)+\bar \alpha \cdot\lambda+\tfrac12\lambda^2\quad (\lambda\geq 0).
\]
Apply \cref{prop:PostComp}(a) to obtain the desired result. 

Part (b). Apply \cref{prop:PostComp}(b), observing that $P_{\epi f}(\bar x,\bar \alpha)=P_1 \delta_{\epi f}(\bar x,\bar \alpha)$ and  $P_1(\lambda\cdot\psi)(\bar x,\bar \alpha)=[P_1(\lambda  f)(\bar x),\bar \alpha+\lambda]$ for all $\lambda\geq 0$~\cite[Theorem~6.6]{Bec 17}. The  fact that  $\bar \lambda>0$ is due to  \cref{prop:PostComp}(c).

 \hfill\Halmos\endproof
\begin{remark}[Prior work] 
  The level-set projection result \Cref{cor:LevelSet} encompasses the result described by Beck~\cite[Theorem~6.30]{Bec 17}. For epigraphical projection, \cref{cor:Epi} generalizes Beck~\cite[Theorem~6.36]{Bec 17} to include functions that aren't finite-valued.
  For functions $f\in \Gamma_0(\bE)$ with open domain, Chierchia et al.~\cite[Proposition~1]{CPP 15} describe an alternative formula for epigraphical projections via proximal maps.
\end{remark}

 \subsection{An SC\boldmath{$^1$}\! optimization framework}\label{sec:PC1}

In this section we present a unified algorithmic framework  for  computing projections onto the level sets and the epigraph of a closed proper convex function. \Cref{cor:LevelSet,cor:Epi}, respectively, guide us in how to compute these projections. For a given $f\in \Gamma_0(\bE)$ and $(\bar x,\bar \alpha)\in \bE\times \R$ such that $f(\bar x)>\bar \alpha$, the epigraphical and level-set projections, respectively, correspond to the proximal map of $f$ with parameter $\lambda$ that solves the scalar problem
\begin{equation}\label{eq:min-theta}
\min_{\lambda\geq 0}\ \theta_\xi(\lambda)\quad (\xi\in \{\epi, \lev{}{}\}),
\end{equation}
for  $\theta_\xi:\R\to\rbar$ given by 
\begin{equation}\label{eq:thetaLev}
\theta_{\xi}(\lambda)=
\begin{cases}
	\bar \phi^f_{\bar x}(\lambda) + \bar \alpha \lambda & \mbox{if $\xi=\lev{}{}$,}\\
    \bar \phi^f_{\bar x}(\lambda) + \bar \alpha \lambda+ \tfrac{1}{2}\lambda^2& \mbox{if $\xi=\epi$.}
\end{cases}
\end{equation}
\Cref{cor:Phi} asserts that $\theta_\xi$ is convex, continuous (possibly in an extended real-valued sense), and continuously differentiable with monotonically increasing, locally Lipschitz 
derivative on  $\R_{++}$. In particular, for any $\lambda>0$,
\begin{equation} \label{eq:theta-prime}
\theta_\xi'(\lambda)=\begin{cases}
	-\eta_{\bar x}^f(\lambda) + \bar \alpha & \mbox{if $\xi=\lev{}{}$,}\\
    -\eta_{\bar x}^f(\lambda) + \bar \alpha +\lambda& \mbox{if $\xi =\epi$,}
\end{cases} 
\end{equation}

The minimization of $\phi_\eta$ could be accomplished using bisection if an upper bound on the optimal $\lambda$ is available. However, the semismoothness of the derivative~\eqref{eq:theta-prime}, described by \cref{prop:SemiEta}, allows us to tap into the powerful SC$^1$ optimization framework~\cite{FaP 03, PaQ 95} that operates on functions $\theta:\R\to\rbar$ that are {\em semismoothly differentiable} (i.e., SC$^1$), which means that at points $\bar \lambda\in\inter(\dom \theta)$, the gradient $\theta'$ exists, and it is locally Lipschitz around $\bar\lambda$ and  semismooth at $\bar \lambda$. The semismooth method, outlined by \cref{alg:General}, applies to the problem~\eqref{eq:min-theta} whenever conditions (A1) and (A2) of Pang and Qi~\cite{PaQ 95} hold, which is the case when $\bar x\in\dom \p f$; see \cref{cor:Phi}.  
 
\Cref{alg:General} uses the notion of a \emph{Bouligand subdifferential}, which for a function $\phi:\R^n\to \rbar$ that is locally Lipschitz at a point $\bar x \in\inter (\dom \phi)$, is defined at $\bar x$ as
$
  \p_B \phi(\bar x)=\set{v}{\exists \{x_k\in D_\phi,\ x_k\to\bar x\}: \nabla \phi(x_k)\to v},
$ 
where $D_\phi$ is the set of points of differentiability of $\phi$. The  {\em Clarke subdifferential} \cite{Cla 83} of $\phi$ at $\bar x$ is 
$
\p_C\phi(\bar x):=\co  \p_B \phi(\bar x),
$
which coincides (on the interior of $\dom \phi$) with the convex subdifferential if $\phi$ is convex.
 

\begin{algorithm}[t]
\begin{algorithmic}
\item[(S.0)] Choose $\lambda_0,\delta>0$, $\{\varepsilon_k\}\downarrow 0$,  and let $\beta,\sigma\in  (0,1)$. Set $k:=0$.
\item[(S.1)] If $|\theta'(\lambda_k)|\leq \delta:$ STOP. 
\item[(S.2)] Choose $g_k\in \p_B(\theta_\xi')(\lambda_k)$ and set 
\[
\Delta_k:=P_{[-\lambda_k,\infty)}\left( -\frac{\theta_\xi'(\lambda_k)}{g_k+\varepsilon_k}\right).
\]
\item[(S.3)] Set 
\[
t_k:=\max_{l\in \bN_0}\set{\beta^l}{\theta_\xi(\lambda_k+\beta^l\Delta_k)\leq \theta_\xi(\lambda_k)+\beta^l\sigma \theta_\xi'(\lambda)\Delta_k}.
\]
\item[(S.4)] Set $\lambda_{k+1}:=\lambda_k+t_k\Delta_k$, $k\leftarrow k+1$, and go to (S.1).
\end{algorithmic}
\caption{SC$^1$ Newton method for minimizing $\theta_\xi$}
\label{alg:General}
\end{algorithm}

\begin{remark}\label{rem:Clarke} Because $\theta_{\xi}$ is convex  and differentiable with  locally Lipschitz derivative on $\R_{++}$, all elements in the  Clarke subdifferential $\p_C(\theta_\xi')(\lambda)$ are nonnegative for all $\lambda>0$~\cite{FaP 03}. In the epigraphical case (i.e., $\xi=\epi$),  the quadratic term in the expression for $\theta_\epi$ in~\eqref{eq:thetaLev} implies that the elements are bounded below by $1$. Thus, the sequence of regularization parameters $\{\varepsilon_k\}\downarrow 0$ in \cref{alg:General} is not necessary, and in fact, if $\theta'$ is piecewise affine, the regularization could be eliminated by setting the constant regularization $\varepsilon_k:=0$ for all $k$, which would improve numerical convergence regardless of the optimality parameter $\delta>0$.
\end{remark}

\subsubsection{The case where \boldmath{$\theta_\xi'$} is concave on \boldmath{$(0,\lambda_l)$}} 

\Cref{cor:LevelSet,cor:Epi} imply that there exists positive parameters $\lambda_l\leq\lambda_u$ such that
\begin{equation} \label{eq:SolLam}
[\lambda_l,\lambda_u]=\argmin_{\lambda\geq 0}\theta_{\xi}=\set{\lambda>0}{\theta_\xi'(\lambda)=0},
\end{equation}
for both the epigraphical and level-set cases. In the epigraphical case in particular, the solution is unique, and thus $\lambda_u=\lambda_l$; see \cref{cor:Epi}(b). If the derivative $\phi_\xi$ is concave on the interval $(0,\lambda_\ell)$, it is possible to take a full Newton step at every iteration while respecting positivity of the iterates, thus saving the computational cost of a backtracking line-search. The simplified iteration is described by \cref{alg:fullNewt}.

For many important functions, e.g., the 1-norm or negative $\log$,  (and their spectral counterparts), the respective map $\theta_\xi'$ is concave on  $\R_{++}$, but, as suggested above,  we only need the following: 

\begin{assumption}[Concavity \bm{$(0,\lambda_l)$}]\label{ass:concave}
	The function $\theta'_\xi$ is concave on $(0,\lambda_l)$.
\end{assumption}

\begin{algorithm}[t]
\begin{algorithmic}
\item[(S.0)] Choose $\lambda_0 >0$, $\delta>0$, and $\{\varepsilon_k\} \downarrow 0$.  Set $k:=0$.
\item[(S.1)] If $|\theta'_\xi(\lambda_k)|\leq \delta:$ STOP. 
\item[(S.2)] Choose $g_k\in \p_{C}(\theta'_\xi)(\lambda_k)$ and set 
\[
\Delta_k:=\max\left\{\frac{-\lambda_k}{2},\ \frac{-\theta'_\xi(\lambda_k)}{g_k+\varepsilon_k}\right\}.
\]
\item[(S.3)] Set $\lambda_{k+1} := \lambda_k + \Delta_k$, $k\leftarrow k+1$, and go to (S.1).
\end{algorithmic}
\caption{Full-step SC$^1$ Newton method}
\label{alg:fullNewt}
\end{algorithm}

\begin{proposition}[Convergence of \cref{alg:fullNewt}]\label{prop:convexCase}
  Under \cref{ass:concave}, the full-step Newton method from \cref{alg:fullNewt} converges to a minimizer of $\theta_\xi$.
\end{proposition}

\proof{Proof.}
 Set $\theta=\theta_\xi$. If $0<\lambda_k<\lambda_l$ for some $k\in \bN$, then by \cref{cor:eta}(a), $\theta'(\lambda_k)<0$ by monotonicity of $-\theta'$. Therefore,
\[
\lambda_{k+1}=\lambda_k-\frac{\theta'(\lambda_k)}{g_k+\varepsilon_k}>\lambda_k.
\]
Since $ -(g_k+\varepsilon_k)$ is a convex subgradient of $-(\theta'+\varepsilon_k(\cdot))$, the concavity of $\theta'_\xi$ implies that
\[
-\theta'(\lambda_{k+1})-\varepsilon_k(\lambda_{k+1}-\lambda_k)\geq -\theta'(\lambda_k)-(\lambda_{k+1}-\lambda_k)(g_k+\varepsilon_k)=0,
\]
and hence $\theta'(\lambda_{k+1})<0$, thus $0<\lambda_k<\lambda_{k+1}<\lambda_l$. Consequently, by an inductive argument, $\{\lambda_k\}$ converges to some $\tilde\lambda$. Therefore,  the sequence $\{g_k\in \p_C(\theta_\xi)(\lambda_k)\}$ is bounded, and hence
\begin{equation*}
0 = (\lambda_{k+1} - \lambda_{k})(g_k + \varepsilon_k) + \theta'(\lambda_k) \to
 \theta'(\tilde{\lambda}),
\end{equation*}
which shows that $\tilde \lambda $ has the desired properties. We hence still  need to cover the case where 
$ \lambda_l<\lambda_k$ for all $k\in\bN$.   In view of~\eqref{eq:SolLam}, we can assume that $\lambda_u<\lambda_k$ for all $k\in\bN$. (Otherwise, a solution has already been obtained.) Since $\theta'(\lambda_k)> 0$ here,  we observe that
\[
0<\lambda_u<\lambda_{k+1} =\lambda_k+ \max\left\{\frac{-\lambda_k}{2},\ \frac{-\theta'(\lambda_k)}{g_k+\varepsilon_k}\right\}\leq \lambda_k,
\]
hence the sequence $\{\lambda_k\}$ converges  to some $\hat \lambda$.   In particular, $\lambda_{k+1}=\tfrac12\lambda_k$ only finitely many times. Hence, without loss of generality, 
$
0= (\lambda_{k+1} - \lambda_{k})(g_k + \varepsilon_k) + \theta'(\lambda_k)\to \theta'(\hat \lambda),
$
which gives $\theta'(\hat \lambda)=0$ also here.
 \hfill\Halmos\endproof



The next example illustrates that cycling may occur in \cref{alg:fullNewt} if Assumption \ref{ass:concave} fails.

\begin{example}[Cycling]\label{ex:Cycle} 
Consider the scalar function $f(x) = 2|x| + \delta_{\left[-1,1\right]}(x)$, and the task of projecting the $(\bar{x},\bar{\alpha}) = (4,-1)$ onto $\epi f$. \Cref{fig:cycle-example} illustrates the function $\theta'_\epi$ whose root we seek.  Then for $\lambda_0$ outside of the interval $\left[1.5,2\right]$ the iterates $\lambda_k\;(k\in\bN)$ generated by \cref{alg:fullNewt} oscillate between $1.5$ and $3$.
\end{example}

\subsection{Numerical Experiments}\label{sec:numerics}

We present numerical experiments that hint at the computational effectiveness of the SC$^1$ optimization framework described in \cref{sec:PC1}. The two experiments in this section were run on an Apple Macbook Air with a 1.8GHz Intel Core i5 and 8Gb RAM running OS 10.14.6. The code was written in C and available at \url{https://github.com/arielgoodwin/epi-proj}.

\subsubsection{Level-set projection: 1-norm} An important instance of the level-set case ($\xi=\lev{}{}$) is the projection onto the unit 1-norm ball $\lev{1}{\|\cdot\|_1}=\set{x\in\R^n}{\|x\|_1\leq 1}$. The derivative of the corresponding function $\theta_{\lev{}{}}$ reads
\[
  \theta_{\lev{}{}}'(\lambda)
  =
  \begin{cases}
    1 - \sum_{i = 1}^n\max\{|x_i| - \lambda, 0\} & \mbox{if $\lambda \geq 0$,}\\
    1-||x||_1 & \mbox{if $\lambda < 0$,}
\end{cases}
\]
which is concave on $\R_+$ (as required) and piecewise affine, as shown by \cref{fig:l1Proj}.
\begin{figure}
	\centering
	\tikzset{>=latex}
	\begin{minipage}[b]{.45\textwidth}
		\centering\small
		\begin{tikzpicture}[scale = 0.9,
			declare function={
			func(\x) = (\x <= 0) * (-6.1) +
			and(\x > 0, \x <= 0.8) * (4*\x - 6.1) +
			and(\x > 0.8, \x <= 1.3) * (3*\x -5.3) +
			and(\x > 1.3, \x <= 2) * (2*\x -4) +
			and(\x > 2, \x <= 3) * (\x - 2) +
			(\x > 3) * (1)
			;
			}
			]
			\begin{axis}[
				axis x line=middle, axis y line=middle,
				ymin=-6.5, ymax=2, ytick={-5,-3,-1,1}, ylabel=$$,
				xmin=-4, xmax=5, xtick={-5,-3,-1,1,3,5}, xlabel=$$,
				domain=-3:5,samples=101, 
				]
				\addplot[thick] {func(x)};
			\end{axis}
		\end{tikzpicture}
		\captionof{figure}{The function $\theta'_\epi$ corresponding to the projection of point $\bar x = (-2,0.8,3,1.3)$ onto the 1-norm unit ball.}
		\label{fig:l1Proj}
	\end{minipage}%
	\hfil
	\begin{minipage}[b]{.45\textwidth}
		\centering\small
	    \tikzset{>=latex}
		\begin{tikzpicture}[scale = 0.9]
			\draw[thin,->] (-1,0) -- (3.5,0) ;
			\draw[thin,->] (0,-2.5) -- (0,2.5);
			\draw[very thick](0.4,-2.6) -- (1.5,-1.5) -- (2,1) -- (3,2) node[anchor=west] {$\phi'_\xi$};
			\draw[dotted,very thick] (2,1) --  (1,0);
			\draw[dotted,very thick] (0.4,-2.6) -- (3,0);
			\foreach \x in {1,2,3}
			\draw (\x cm,-1pt) -- (\x cm,1pt) node[anchor=north] {$\x$};
			\foreach \y in {-2,-1,1,2}
			\draw (-1pt,\y cm) -- (1pt,\y cm) node[anchor=east] {$\y$};
		\end{tikzpicture}
		\captionof{figure}{The function $\theta'_\epi$ for \cref{ex:Cycle}, for which \cref{alg:fullNewt} may cycle.\\}
		\label{fig:cycle-example}
	\end{minipage}
\end{figure}

We implemented \cref{alg:fullNewt} and compared it numerically to two state-of-the-art algorithms specifically tailored to 1-norm-ball projection, namely Condat's sorting-based method~\cite{Con 16} as implemented in the code \texttt{condat\_l1ballproject.c}, and Liu and Ye's improved bisection algorithm (IBIS)~\cite{Liu 09} implemented in the \texttt{eplb} module in SLEP~\cite{SLEP}.

\begin{table}[t]
\centering
\begin{tabular}{cccc|ccc}
	 \toprule[1pt]
    {$n$} 		& {Algorithm 2} 			& {Condat} 	    & {IBIS}				&  {Algorithm 2} 			& {Condat}  	  	& {IBIS}			
    \\ \midrule[.5pt]
	& \multicolumn{3}{c|}{$\sigma = 0.1$} 					& \multicolumn{3}{c}{$\sigma = 0.05$}
    \\ \midrule[.5pt]
    {$20$} 		&\num{1.94e-6}		& \num{1.53e-6} 	&\num{1.83e-6}		& \num{1.93e-6}	& \num{1.41e-6} 	&\num{1.99e-6}		\\
    {$10^3$} 	&\num{3.33e-5}		& \num{2.11e-5} 	&\num{3.65e-5}		& \num{3.38e-5}	& \num{2.23e-5}	&\num{4.15e-5}		\\
	{$10^6$} 	&\num{2.08e-2}		& \num{1.44e-2} 	&\num{2.89e-2}		& \num{2.18e-2} 	&\num{1.44e-2}		&\num{3.42e-2}		\\[6pt]
	\midrule[.5pt]
    & \multicolumn{3}{c|}{$\sigma = 0.01$}	
    							& \multicolumn{3}{c}{$\sigma = 0.005$} 				\\
   \midrule[.5pt]							
 {$20$} 	    &\num{2.05e-6}		&\num{1.45e-6}		&\num{1.87e-6}		& \num{1.92e-6}	& \num{1.36e-6} 	&\num{2.32e-6}		\\
 {$10^3$}     &\num{3.14e-5}		&\num{2.57e-5}		&\num{4.07e-5}		& \num{3.06e-5}	& \num{2.68e-5}	&\num{4.46e-5}		\\
 {$10^6$}     &\num{1.93e-2}		&\num{1.48e-2} 	&\num{3.73e-2}		& \num{1.89e-2} 	&\num{1.50e-2}		&\num{4.00e-2}		\\ \bottomrule[1pt]
\end{tabular}
\medskip
\caption{Average time (seconds) for projecting vectors onto the 1-norm unit ball in dimension $n$, with coordinates chosen using Gaussian distributions with standard deviation $\sigma$.\label{Tab:numexp}}
\end{table}

The entries of the projected vectors $\bar x\in\R^n$ are drawn from a Gaussian distribution with zero mean and standard deviations $\sigma=\{0.1,0.05,0.01,0.005\}$. The optimality tolerance was fixed at $\delta=10^{-15}$, as in step (S.1) of \cref{alg:fullNewt}.
\Cref{Tab:numexp} reports the average time required to compute the projection over $10^5$ trials for vectors of dimension $n \in \{20, 10^3\}$, and over 500 trials for $n=10^6$. The initial point $\lambda_0>0$ \cref{alg:fullNewt} was chosen by sampling $\sqrt{n}\log n$ coordinates randomly from the vector $\bar x$ and setting $\lambda_0$ to be the largest of their absolute values. 
Observe that \cref{alg:fullNewt} exhibits comparable performance relative to the specialized algorithms.

\subsubsection{Level-set projection: negative sum-log}

\begin{figure}[t]
\centering
\begin{tabular}{@{}cc@{}}
	\begin{tikzpicture}[scale = 0.85,
		declare function={
		func(\x) = (\x <= 0) * (\x - 1) +
				(\x > 0) * (\x - 1 + ln((1+sqrt(1+4*\x))/2))
		;
		}
	]
	\begin{axis}[
  	axis x line=middle, axis y line=middle,
  	ymin=-4, ymax=6, ytick={-3,-1,...,5}, ylabel={},
  	xmin=-3, xmax=5, xtick={-3,-1,...,5}, xlabel={},
  	domain=-3:5,samples = 100, 
	]
	\addplot [thick] {func(x)};
	\end{axis}
	\end{tikzpicture}
&
	\begin{tikzpicture}[scale = 0.85,
		declare function={
		func(\x) = (\x > 0) * (\x - 1 + ln((-1+sqrt(1+4*\x))/2)) +
				(\x <= 0) * (-10)
		;
		}
	]
	\begin{axis}[
  	axis x line=middle, axis y line=middle,
  	ymin=-4, ymax=6, ytick={-3,-1,...,5}, ylabel={},
  	xmin=-2, xmax=5, xtick={-1,0,...,5}, xlabel={},
  	domain=-1:5,samples=100, 
	]
	\addplot [thick] {func(x)};
	\end{axis}
	\end{tikzpicture}
	\\[6pt] $(\bar x, \bar \alpha)=(+1,-1)$ & $(\bar x, \bar \alpha)=(-1,-1)$
\end{tabular}
	\caption{The graph of the function $\theta'_\epi(\lambda)$ that corresponds to the base points $(\bar x,\bar \alpha)$ shown for each figure. The left panel depicts the case where $P_{\cl(\dom f)}(\bar x)\in\dom f$; the right panel depicts the case where $P_{\cl(\dom f)}(\bar x)\notin\dom f$. \label{fig:closdom}}
\end{figure}
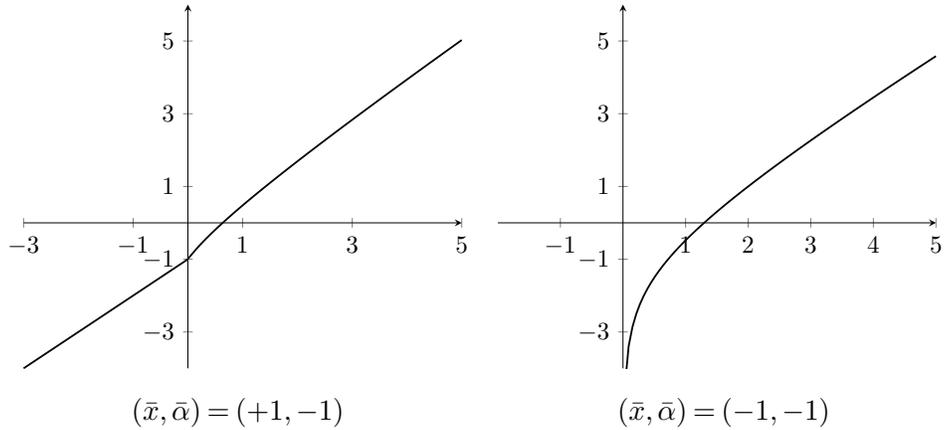

We now consider the epigraphical projection for a function that is not polyhedral. Define the function $f:x\in\R^n\mapsto - \sum_{i = 1}^n\log x_i$, where we take the negative logarithm to be $\infty$ outside the positive orthant. \Cref{fig:closdom} illustrate the function $\phi'_\epi$ for the case when $P_{\cl(\dom f)}(\bar x)$ is in, and not in, the domain of $f$. These functions are concave over $(0,\infty)$. Hence $-\theta'_\xi$ is convex over this interval and \cref{alg:fullNewt} applies.

We numerically compare \cref{alg:fullNewt} and the bisection method as solution approaches for~\eqref{eq:min-theta}. The coordinates of $\bar x$ were chosen uniformly at random on the interval $[-1,1]$, and the value $\bar \alpha$ was chosen uniformly at random on the interval $[-2,-0.5]$. The initial value $\lambda_0$ was chosen to be $\sqrt N$. The termination condition for \cref{alg:fullNewt} was $|\theta'_\xi(\lambda)| < 10^{-4}$, and the termination conditions for bisection was $|\theta'_\xi(\lambda)| < 10^{-4}$ (labeled \emph{Bisection 1}) and $|b-a| < 10^{-8}$ (labeled \emph{Bisection 2}), where $a,b$ denote the endpoints of the bisection interval. \Cref{tab:log} shows the average times over $10^5$ trials when $n \in \{1, 10^3\}$, and over 500 trials when $n=10^6$.

\begin{table}
\centering
\begin{tabular}{lccc}
	\toprule[1pt]
					& $n = 1$			& $n = 10^3$		& $n = 10^6$ \\ \midrule[.5pt]
{SSN}				& \num{8.76e-7}	&\num{1.69e-4}		&\num{1.90e-1}		      \\
{Bisection 1}         	& \num{2.36e-6}       	&\num{1.88e-3}		&\num{2.89e0}                  \\ 
{Bisection 2}		& \num{2.66e-6}        &\num{1.08e-3}	&\num{1.16e0} \\ \bottomrule[1pt]
\end{tabular}
\caption{Time (seconds) for projecting vectors onto the epigraph of $f(x)=-\sum_{i=1}^n\log x_i$ in various dimensions $n$.\label{tab:log}}
\end{table}

\subsubsection{Discussion}\label{sec:Discuss}

The numerical examples we presented extend easily to other useful cases involving matrices, such as the nuclear norm on $\R^{m\times n}$ and the barrier function $-\log\det $ on the space of symmetric matrices, using variational formulas that depend on matrix spectra~\cite{Lew 95, Lew 96}.  In these cases, the main computational effort involves computing singular value and eigenvalue decompositions, respectively, of the matrix iterates. 

The cases where $\theta'_\xi$ does not satisfy either Assumption \ref{ass:concave} or the domain condition $P_{\cl(\dom f)}(\bar x) \in \dom f$  lies  outside the theoretical guarantees presented in this section, though the algorithms we present may still work in practice. In the case where $\dom f\subsetneq\bE$ is open, the formula provided by Chierchia et al.~\cite[Proposition~1]{CPP 15} is a viable option.

\section{Final remarks}

Our analysis on the variational properties of epigraphical projections and infimal convolution is motivated by the authors' larger research interests on variations of first-order methods that operate in a lifted space. The promising work by Chierchia et al.~\cite{CPP 15} on epigraphical-projection methods for minimizing convex functions over $p$-norm constraints shows promise for this algorithmic approach, and we aim to develop methods for more general problem classes. We are also motivated by statistical M-estimation approaches that include as an additional unknown a particular parameter that characterizes data distribution~\cite{CoM 20}. The variational calculus that we derive is a useful tool for developing algorithmic approaches for solving these lifted M-estimation problems. 


There are at least two avenues of future research that extend our analysis in this paper.

\paragraph{$\bm{K}$-epigraphical projections.} A significant generalization of the post-composition operation defined in \cref{sec:post-composition} occurs when we allow compositions of the form $f=g\circ H:\bE_1\to\rbar$, where
 \begin{itemize}
 \item $K\subset \bE_1$ a  closed convex cone;
 \item $H:\bE_1\to \bE_2$ $K$-convex, i.e., the \emph{$K$-epigraph} $\set{(X,Y)}{Y-H(x)\in K}$ is convex;
 \item $g\in \Gamma_0(\bE_2)$ \emph{$K$-increasing}, i.e., $g \leq g((\cdot)+v)$ for all $v\in K$.
 \end{itemize}
This {\em convex convex-composite} setting was studied by Burke et al.~\cite{BHN 20}, and  the required subdifferential formulas for the analysis are readily available. This may lead to a proximal calculus and ultimately to formulas and algorithms for projecting onto $K$-epigraphs, thus encompassing the study in \cref{sec:post-composition}.

\paragraph{Semismoothness* of subdifferential operators.} The notion of semismooth* sets and maps is recent and still in development. One of the critical conditions in our study is the semismoothness* of the subdifferential operator $\p f$, which also occurs in a recent report by Khanh et al.~\cite{KMP 20}. This suggests an important avenue of research that relaxes the overarching convexity assumption and, in particular, establishes verifiable  sufficient conditions.

\section*{Acknowledgments}   M.P. Friedlander and T. Hoheisel are supported by NSERC Discovery grants, while A. Goodwin's work was partially supported by an NSERC summer research stipend.
T. Hoheisel would like to thank Dr.~Matus Benko,  University of Vienna, for valuable discussions on semismoothness*.


%


%
%
%


\end{document}